\documentclass[reqno,a4paper,11pt]{amsart}
\usepackage{graphicx}
\usepackage{geometry}
\usepackage{a4wide}

\usepackage{latexsym}
\usepackage{amsmath}
\usepackage{amsfonts}
\usepackage{amssymb}
\usepackage{amsthm}
\usepackage{mathrsfs}

\usepackage{xfrac}
\usepackage{tikz}
\usepackage{mathtools}

\usepackage[pdftex]{hyperref}
\usepackage[utf8x]{inputenc}
\usepackage[english]{babel}

\newcommand{\FD}{D_{fr}}
\newcommand{\triend}{\mbox{\hspace{0.2mm}}\hfill$\triangle$}
\newcommand{\black}{\mbox{\hspace{0.2mm}}\hfill$\blacksquare$}

\newtheorem{theorem}{Theorem}
\newtheorem{proposition}[theorem]{Proposition}
\newtheorem{lemma}[theorem]{Lemma}
\newtheorem{corollary}[theorem]{Corollary}
\newtheorem*{theorem*}{Theorem}

\theoremstyle{definition} 
\newtheorem{definition}[theorem]{Definition}

\theoremstyle{remark} 
\newtheorem{remark}[theorem]{Remark}
\newtheorem{example}[theorem]{Example}
\newtheorem{assumption}{Assumption}

\title{Restriction theorems for $\mu$-(semi)stable framed sheaves}

\author{Francesco Sala}

\address{Department of Mathematics, School of Mathematical and Computer Sciences, Heriot-Watt University, Colin Maclaurin Building, Riccarton, Edinburgh EH14 4AS, United Kingdom and Maxwell Institute for Mathematical Sciences, Edinburgh, United Kingdom}

\email{F.Sala@hw.ac.uk, salafra83@gmail.com}

\begin{document}

\begin{abstract}
 We provide a generalization of Mehta-Ramanathan restriction theorems to framed sheaves: we prove that the restriction of a $\mu$-semistable framed sheaf on a nonsingular projective irreducible variety of dimension $d\geq 2$ to a general hypersurface of sufficiently high degree is again $\mu$-semistable. The same holds for $\mu$-stability under some additional assumptions.
\end{abstract}

\maketitle

\setlength{\parskip}{0.2ex} 

\tableofcontents

\setlength{\parskip}{0.6ex plus 0.3ex minus 0.2ex} 

\section{Introduction}

\thispagestyle{empty}

In \cite{art:donaldson1984}, Donaldson proved that the moduli space of gauge-equivalence classes of framed $SU(r)$-instantons with instanton charge $n$ on $S^4$ is isomorphic to the moduli space of isomorphism classes of vector bundles on $\mathbb{CP}^2$ of rank $r$ and second Chern class $n$ that are trivial along a fixed line $l_{\infty}$, and with a fixed trivialization there. It is an open subset $\mathcal{M}^{reg}(r,n)$ of the moduli space $\mathcal{M}(r,n)$ of \textit{framed sheaves} on $\mathbb{CP}^2$, that is, the moduli space parametrizing pairs $(E,\alpha)$, modulo isomorphism, where $E$ is a torsion free sheaf on $\mathbb{CP}^2$ of rank $r$ and $c_2(E)=n$, locally trivial in a neighbourhood of $l_{\infty}$, and $\alpha\colon E\vert_{l_{\infty}}\stackrel{\sim}{\rightarrow} \mathcal{O}_{l_{\infty}}^{\oplus r}$ is the \textit{framing at infini\-ty}. $\mathcal{M}(r,n)$ is a nonsingular quasi-projective variety of dimension $2rn.$ Moreover, it admits a description in terms of linear data, the so-called \
\textit{ADHM} data (see, for example, \cite[Chapter 2]{book:nakajima1999}). In some sense we can look at $\mathcal{M}(r,n)$ as a partial compactification of $\mathcal{M}^{reg}(r,n).$ There exists another type of partial compactification $\mathcal{M}^{Uh}(r,n)$ of the latter moduli space, called \textit{Uhlenbeck-Donaldson compactification}: using linear data and geometric invariant theory it is possible to construct a projective morphism
\begin{equation}\label{eq:pi}
\pi_r\colon \mathcal{M}(r,n)\rightarrow \mathcal{M}^{Uh}(r,n)=\bigsqcup_{i=0}^{n} \mathcal{M}^{reg}(r,n-i)\times \mathrm{Sym}^i(\mathbb{C}^2)
\end{equation}
such that the restriction to the ``locally free'' part is an isomorphism with its image (see \cite[Chapter 3]{book:nakajima1999}).

The moduli spaces $\mathcal{M}(r,n)$ can be regarded as higher-rank generalizations of Hilbert schemes of $n$-points on the complex affine plane. From this point of view, the previous morphism is a higher-rank generalization of Hilbert-Chow morphism for Hilbert schemes of points on the complex affine plane. 

Because of the relation with moduli space of framed instantons, since Nekrasov's partition function was introduced in \cite{art:nekrasov2003}, the moduli space $\mathcal{M}(r,n)$ has been studied quite intensively (see, e.g., \cite{art:bruzzofucitomoralestanzini2003, book:nakajima1999, art:nakajimayoshioka2005-I, art:nakajimayoshioka2005-II,  phd:carlsson2009}) and the geometry of moduli spaces of framed sheaves on the complex projective plane is quite well known.

In \cite{art:huybrechtslehn1995-I, art:huybrechtslehn1995-II} Huybrechts and Lehn laid the foundations of a systematic theory of framed sheaves on varieties of arbitrary dimension (they used different names to denote the same object like, e.g., \textit{stable pairs}, \textit{framed modules}, etc). Let $X$ be a nonsingular, projective, irreducible variety of dimension $d$ defined over an algebraically closed field $k$ of characteristic zero. A \textit{framed sheaf} is a pair $(E,\alpha)$ where $E$ is a coherent sheaf on $X$ and $\alpha$ is a morphism from $E$ to a fixed coherent sheaf $F$, called \textit{framing sheaf}. They define a generalization of Gieseker semistability (resp.\ $\mu$-semistability) for framed sheaves that depends on a polarization and a rational polynomial $\delta$ of degree $d-1$ with positive leading coefficient (resp.\ a rational number $\delta_1$). 

Fix a \textit{numerical polynomial} $P$ of degree $d$, i.e., a rational polynomial $P(n)$ of degree $d$ such that $P(a)\in\mathbb{Z}$ for any integer $a.$ Let us denote by $\underline{\mathcal{M}}^{ss}_\delta(X;F,P)$ (resp.\ $\underline{\mathcal{M}}^{s}_\delta(X;F,P)$) the contravariant functor from the category of Noetherian $k$-schemes of finite type to the category of sets, that associates to every scheme $T$ the set of isomorphism classes of families of semistable (resp.\ stable) framed sheaves with Hilbert polynomial $P$ parametrized by $T.$ The main result in their papers is the following:
\begin{theorem*}[Huybrechts, Lehn]
There exists a projective scheme $\mathcal{M}^{ss}_\delta(X;F,P)$ that corepresents the functor $\underline{\mathcal{M}}^{ss}_\delta(X;F,P).$ Moreover, there is an open subscheme $\mathcal{M}^{s}_\delta(X;F,P)$ of $\mathcal{M}^{ss}_\delta(X;F,P)$ that represents the functor $\underline{\mathcal{M}}^{s}_\delta(X;F,P)$, i.e., $\mathcal{M}^{s}_\delta(X;F,P)$ is a fine moduli space for stable framed sheaves.
\end{theorem*}
If $X$ is a surface, we can extend the original definition of framed sheaves on $\mathbb{CP}^2$ with framing along a fixed line in the following way: let $F$ be a coherent sheaf on $X$, supported on a big and nef curve $D$, such that $F$ is a Gieseker semistable locally free $\mathcal{O}_{D}$-module. A $(D,F)$\textit{-framed sheaf} is a framed sheaf $(E,\alpha\colon E\rightarrow F)$, where $\ker\alpha$ is torsion free, $E$ is locally free in a neighbourhood of $D$ and $\alpha\vert_D$ is an isomorphism. It is possible to prove that there exists a rational polynomial $\delta$ such that there is an open subset in $\mathcal{M}^{s}_\delta(X;F,P)$ parametrizing $(D,F)$-framed sheaves on $X$ with Hilbert polynomial $P.$ Moreover if the surface $X$ is rational and $D$ is a smooth irreducible big and nef curve of genus zero, the moduli space of $(D,F)$-framed sheaves is a nonsingular quasi-projective variety (see \cite{art:bruzzomarkushevich2011}). It is possible to generalize the definition of $(D,F)$-framed sheaves 
to varieties of arbitrary dimension (see Definition \ref{def:locfree}).

Leaving aside the results on the representability of the moduli functor $\underline{\mathcal{M}}_\delta^{(s)s}(X;F,P)$ discussed, a complete theory of framed sheaves and a study of the geometry of their moduli spaces is missing in the literature. In the present paper we fill one of the gaps of the theory, by providing a generalization of the Mehta-Ramanathan restriction theorems:
\begin{theorem*}
Let $X$ be a nonsingular, projective, irreducible variety of dimension $d$, defined over an algebraically closed field $k$ of characteristic zero, endowed with a very ample line bundle $\mathcal{O}_X(1).$ Let $F$ be a coherent sheaf on $X$ supported on a divisor $\FD.$ Let $\mathcal{E}=(E,\alpha\colon E\rightarrow F)$ be a framed sheaf on $X$ of positive rank with nonzero framing. If $\mathcal{E}$ is $\mu$-semistable with respect to $\delta_1$, then there is a positive integer $a_0$ such that for all $a\geq a_0$ there is a dense open subset $U_a\subset \vert \mathcal{O}_X(a)\vert$ such that for all $D\in U_a$ the divisor $D$ is smooth, meets transversely the divisor $\FD$ and $\mathcal{E}\vert_D$ is $\mu$-semistable with respect to $a\delta_1.$
 
The same statement holds with ``$\mu$-semistable'' replaced by ``$\mu$-stable'' under the following additional assumptions: the framing sheaf $F$ is a locally free $\mathcal{O}_{\FD}$-module and $\mathcal{E}$ is a $(\FD,F)$-framed sheaf on $X.$
\end{theorem*}
Mehta-Ramanathan restriction theorems are very useful as they often allow one to reduce a problem from a higher-dimensional variety to a surface or even to a curve, as for example happens with the proof of Hitchin-Kobayashi correspondence (see \cite[Chapter VI]{book:kobayashi1987}). 

The classical Mehta-Ramanathan restriction theorems are also used in the algebro-geometric construction of the Uhlenbeck-Donaldson compactification of moduli space of $\mu$-stable vector bundles on a nonsingular projective surface (see \cite{art:li1997} and \cite[Section 8.2]{book:huybrechtslehn2010}). In the same way, our framed version of these theorems is used in a work of Bruzzo, Markushevich and Tikhomirov in the construction of the Uhlenbeck-Donaldson compactification for framed sheaves (see \cite{art:bruzzomarkushevichtikhomirov2010}). In this way they provide a generalization of the morphism $\pi_r$ (see formula \eqref{eq:pi}) to an arbitrary smooth projective complex surface.

The main difficulty in the generalization of Mehta-Ramanathan restriction theorems has been the lack of some basic tools in the theory of framed sheaves. In this paper, we provide new tools for the study of the (semi)stability condition for framed sheaves. We construct the (relative) Harder-Narasimhan filtration, used in the proof of the first restriction theorem. We also define the Jordan-H\"older filtration and construct the (extended) socle, used in the proof of the second restriction theorem. We want to emphasize that the second theorem is not proved in the same generality as the first one. This is due to some technical problems: for example, in general it is impossible to define a natural framing on the double dual of the underlying sheaf of a framed sheaf; moreover, under the assumption that $F$ is an arbitrary coherent sheaf, if a framed sheaf is simple, it is no more true that it is remains simple upon the restriction to a (general) divisor. To circumvent these difficulties, we had to strengthen the hypotheses.

We follow rather closely the approach chosen by Huybrechts and Lehn in their book \cite{book:huybrechtslehn2010} to study Gieseker and $\mu$-semistability conditions for pure sheaves. As they do in \cite{book:huybrechtslehn2010}, we study the (relative) Harder-Narasimhan filtration and the Jordan-H\"older filtration for the polynomial (semi)stability condition (see from Section 3 to Section 8). A careful check shows that similar results hold also for the $\mu$-semistability condition (cf. Remark \ref{rem:results} and the subsequent theorem). Moreover, unlike \cite{art:huybrechtslehn1995-I, art:huybrechtslehn1995-II}, where a theory of framed sheaves for a smooth projective irreducible variety over an algebraically closed field of characteristic zero is developed, the ambient space is for us a projective scheme over such a field, unless otherwise stated. To introduce suitable (semi)stability conditions depending on Hilbert polynomials, only the projectivity condition is needed.

Each section of the paper starts with a summary which describes when the results in the framed case coincide with the corresponding ones in the nonframed case or when there are unexpected phenomena. We refer to \cite{book:huybrechtslehn2010} for the nonframed case.

\subsection*{Conventions}

A sheaf of $\mathcal{O}_Y$-modules on a Noetherian scheme $Y$ is always meant to be \emph{coherent}, so we shall omit the adjective. As usual in the literature, we identify a \emph{vector bundle} on $Y$ with the sheaf of its sections.

Let $E$ be a sheaf on $Y.$ The \emph{support of} $E$ is the closed set $\mathrm{Supp}(E) := \{x \in Y\,\vert\,E_x \neq 0\}.$ Its dimension is called the \emph{dimension of} $E$ and is denoted by $\dim(E).$ The sheaf $E$ is \emph{pure} if for all nontrivial subsheaves $E'\subset E$, we have $\dim(E')=\dim(E).$ Recall that the \emph{torsion filtration} of $E$ is the unique filtration
\begin{equation*}
0\subset T_0(E)\subset \cdots\subset T_{\dim(E)-1}(E)\subset T_{\dim(E)}(E)=E,
\end{equation*}
where $T_i(E)$ is the maximal subsheaf of $E$ of dimension less or equal than $i$, for $i=0, \ldots, \dim(E).$ Thus $E$ is pure if and only if $T_{\dim(E)-1}(E)=0.$

If $Y$ is integral, a sheaf $E$ on $Y$ is said to be \emph{torsion free} if for any $y\in Y$ and $s\in\mathcal{O}_{Y,y}\setminus\{0\}$, the multiplication morphism $s\colon E_y\rightarrow E_y$ by $s$ is injective. This condition is equivalent to $T_{\dim(Y)-1}(E)=0.$

Let $Y$ be a projective scheme over a field. Recall that the Euler characteristic of a sheaf $E$ is $\chi(E):=\sum_i (−1)^i \dim \mathrm{H}^i(Y,E).$ Fix an ample line bundle $\mathcal{O}(1)$ on $Y.$ Let $P(E,n):=\chi(E\otimes \mathcal{O}(n))$ be the \emph{Hilbert polynomial of} $E.$ By \cite[Lemma 1.2.1]{book:huybrechtslehn2010}, $P(E,n)$ can be uniquely written in the form
\begin{equation*}
P(E,n)=\sum_{i=0}^{\dim(E)}\beta_i(E)\frac{n^i}{i!},
\end{equation*}
where $\beta_i(E)$ are rational coefficients. Moreover for $E\neq 0$, the leading coefficient $\beta_{\dim(E)}(E)$, called \emph{multiplicity} of $E$, is positive. Note that $\beta_{\dim(Y)}(\mathcal{O}_Y)$ is the degree $\deg(Y)$ of $Y$ with respect to $\mathcal{O}(1).$ The \emph{reduced Hilbert polynomial} $p(E,n)$ of $E$ is
\begin{equation*}
p(E,n):=\frac{P(E,n)}{\beta_{\dim(E)}(E)}. 
\end{equation*}
We call \emph{hat-slope} the quantity 
\begin{equation*}
\hat{\mu}(E)=\frac{\beta_{\dim(E)-1}(E)}{\beta_{\dim(E)}(E)}.
\end{equation*}
For a polynomial $\displaystyle P(n)=\sum_{i=0}^t \beta_i(n^i/i!)$ with $\beta_t\neq 0$, we define $\hat{\mu}(P)=\beta_{t-1}/\beta_t.$

For a sheaf $E$ on $Y$, we call \emph{rank} of $E$ the quantity
\begin{equation*}
\mathrm{rk}(E):=\frac{\beta_{\dim(Y)}(E)}{\deg(Y)}. 
\end{equation*}
In general, $\mathrm{rk}(E)$ is not an integer, but if $Y$ is integral, $\mathrm{rk}(E)\in\mathbb{Z}_{\geq 0}.$ The \emph{degree} of $E$ is
\begin{equation*}
\deg(E):=\beta_{\dim(Y)-1}(E)-\mathrm{rk}(E)\beta_{\dim(Y)-1}(\mathcal{O}_Y).
\end{equation*}
Note that for a sheaf $E$ of dimension $\dim(Y)-1$, $\deg(E)=\beta_{\dim(Y)-1}(E)$, for any sheaf of dimension less than $\dim(Y)-1$, its degree is zero.

If $E$ is a sheaf of dimension $\dim(Y)$, its \emph{slope} is
\begin{equation*}
\mu(E):=\frac{\deg(E)}{\mathrm{rk}(E)}. 
\end{equation*}
Assume that $Y$ is a nonsingular projective irreducible variety and let $E$ be a sheaf on it. By the Hirzebruch-Riemann-Roch theorem, the \emph{degree} $\deg(E)$ of $E$, introduced before, coincides with $c_1(E)\cdot H^{\dim(Y)-1}$, where $H\in \vert \mathcal{O}(1)\vert$ is a hyperplane section. In particular, $\deg(E)=\deg(\det(E))$, where $\det(E)$ is the \emph{determinant line bundle} of $E$ (cf. \cite[Section 1.1.17]{book:huybrechtslehn2010}). 

Let us denote by $k$ an algebraically closed field of characteristic zero. A \emph{polarized scheme} is a pair $(X,\mathcal{O}_X(1))$, where $X$ is a projective scheme, defined over $k$, and $\mathcal{O}_X(1)$ a very ample line bundle on it. A \emph{nonsingular polarized variety} is a polarized scheme $(X,\mathcal{O}_X(1))$, where $X$ is a nonsingular irreducible variety.

Let $Y\rightarrow S$ be a morphism of finite type of Noetherian schemes. If $T\rightarrow S$ is an $S$-scheme, we denote by $Y_T$ the fibre product $T\times_S Y$ and by $p_T\colon Y_T\rightarrow T$ and $p_Y\colon Y_T\rightarrow Y$ the natural projections. If $E$ is a sheaf on $Y$, we denote by $E_T$ its pull-back to $Y_T.$ We use similar notation for morphisms between sheaves on $Y.$

For $s\in S$ we denote by $Y_s$ the fibre $\mathrm{Spec}(k(s))\times_S Y$. For a sheaf $E$ on $Y$, we denote by $E_s$ its pull-back to $Y_s.$ Often, we shall think of $E$ as a collection of sheaves $E_s$ parametrized by $s \in S.$ If $\alpha\colon E\rightarrow F$ is a morphism of sheaves on $Y$, $\alpha_s$ denotes its pull-back to $Y_s.$

Let $g\colon Y\rightarrow S$ be a projective morphism. A $g$\emph{-ample} line bundle is a line bundle on $Y$ such that the restriction to any fibre $Y_s$ is ample for any $s\in S.$

Let $S$ be an integral $k$-scheme of finite type, $f\colon \mathcal{X}\rightarrow S$ a projective morphism with equidimensional fibres and $\mathcal{O}_\mathcal{X}(1)$ an $f$-ample line bundle. It follows that for any $s\in S$ the pair $(\mathcal{X}_s,\mathcal{O}_\mathcal{X}(1)_s)$ is a polarized scheme of a fixed dimension $d.$ We shall call the pair $(f\colon \mathcal{X}\rightarrow S, \mathcal{O}_\mathcal{X}(1))$ \emph{relative polarized scheme of dimension} $d.$

\subsection*{Notation}

As usual in the literature, we shall use $\square$ to close the proofs, except when inside a proof of a statement, we need to close the proof of a intermediate result. In this case we shall use $\blacksquare$ to close the proof of the intermediate statement and $\square$ to close the main proof. This is, for example, the case of the proofs of Proposition \ref{prop:dest} and Theorem \ref{thm:mr1}. Finally, we shall use $\triangle$ to close Remarks and Examples. 

\subsubsection*{Acknowledgements}

This paper is largely based upon the author’s PhD thesis \cite{phd:sala2011}. The author thanks his supervisors Ugo Bruzzo and Dimitri Markushevich for suggesting this problem and for their constant support. Also, he is indebted to the anonymous referee, whose several useful remarks helped him to substantially improve the results of the paper. He also thanks Daniel Huybrechts, Thomas Nevins and Edoardo Sernesi for useful suggestions and Luis \'{A}lvarez-C\'{o}nsul for interesting discussions. This paper was mostly written while the author was staying at SISSA and Université Lille 1. The last draft of the paper was written while the author was staying at Heriot-Watt University in Edinburgh. He thanks those institutions for hospitality and support. He was partially supported by PRIN ``Geometria delle varietà algebriche e dei loro spazi di moduli'', co-funded by MIUR (cofin 2008), the grant of the French Agence Nationale de Recherche VHSMOD-2009 Nr. ANR-09-BLAN-0104, the European Research Network ``GREFI-GRIFGA'' and the ERASMUS ``Student Mobility for Placements'' Programme.

\section{Preliminaries on framed sheaves}

In this section we introduce the notion of \textit{framed sheaf} and \textit{morphism of framed sheaves}. Moreover for such objects we introduce some invariants, like the \textit{framed Hilbert polynomial} and the \textit{framed degree}. When the framing is zero, a framed sheaf is just its underlying sheaf and these notions coincide with the classical ones (see \cite[Section 1.2]{book:huybrechtslehn2010}).

Let $(X,\mathcal{O}_X(1))$ be a polarized scheme. Fix a positive integer $d\leq \dim(X)$, a sheaf $F$ on $X$ and a polynomial 
\begin{equation*}
\delta(n):=\delta_1 \frac{n^{d-1}}{(d-1)!}+\delta_2 \frac{n^{d-2}}{(d-2)!} +\cdots+\delta_d \in\mathbb{Q}[n]
\end{equation*} 
with $\delta_1>0.$ We call $F$ \textit{framing sheaf} and $\delta$ \textit{stability polynomial}.
\begin{definition}
A \textit{framed sheaf} on $X$ is a pair $\mathcal{E}:=(E, \alpha)$, where $E$ is a sheaf on $X$ and $\alpha\colon  E\rightarrow F$ is a morphism of sheaves. We call $\alpha$ \textit{framing} of $E$.
\end{definition}
For any framed sheaf $\mathcal{E}=(E,\alpha)$, we define the function $\epsilon(\alpha)$ by
\begin{equation*}
\epsilon(\alpha) := \left\{ \begin{array}{ll}
         1 & \mbox{if $\alpha\neq 0$},\\
        0 & \mbox{if $\alpha=0$}.\end{array} \right.
\end{equation*}
The \textit{framed Hilbert polynomial of} $\mathcal{E}$ is
\begin{equation*}
P(\mathcal{E},n):= P(E,n)-\epsilon(\alpha)\delta(n),
\end{equation*}
and its \textit{reduced framed Hilbert polynomial} is
\begin{equation*}
p(\mathcal{E},n):=\frac{P(\mathcal{E},n)}{\beta_{\dim(E)}(E)}.
\end{equation*}
We shall call \textit{Hilbert polynomial} (resp.\ \textit{dimension}) of $\mathcal{E}$ the Hilbert polynomial $P(E)$ (resp.\ the dimension $\dim(E)$) of $E.$ 

If $E'$ is a subsheaf of $E$ with quotient $E'':=E/E'$, the framing $\alpha$ induces framings $\alpha':=\alpha\vert_{E'}$ on $E'$ and $\alpha''$ on $E''$, where the framing $\alpha''$ is defined in the following way: $\alpha''=0$ if $\alpha'\neq 0$, else $\alpha''$ is the induced morphism on $E''.$ With this convention the framed Hilbert polynomial of $\mathcal{E}$ behaves additively:
\begin{equation}\label{eq:somma}
P(\mathcal{E})=P(E', \alpha')+P(E'', \alpha'').
\end{equation}
\textbf{Notation:} If $\mathcal{E}=(E,\alpha)$ is a framed sheaf on $X$ and $E'$ is a subsheaf of $E$, then we denote by $\mathcal{E}'$ the framed sheaf $(E', \alpha')$ and by $\mathcal{E}/E'$ the framed sheaf $(E'', \alpha'')$.

Thus we have a canonical framing on subsheaves and on quotients. The same happens for subquotients, indeed we have the following result.
\begin{lemma}{\normalfont (\cite[Lemma 1.12]{art:huybrechtslehn1995-II}).}\label{lem:fra}
Let $H\subset G\subset E$ be sheaves and $\alpha$ a framing of $E$. Then the framings induced on $G/H$ as a quotient of $G$ and as a subsheaf of $E/H$ agree. Moreover
\begin{equation*}
P\Big(\frac{\mathcal{E}/H}{G/H}\Big)=P\left(\mathcal{E}/G\right).
\end{equation*}
\end{lemma}
Now we introduce the notion of a morphism of framed sheaves.
\begin{definition}
Let $\mathcal{E}=(E, \alpha)$ and $\mathcal{G}=(G, \beta)$ be framed sheaves. A \textit{morphism of framed sheaves} $\varphi\colon \mathcal{E}\rightarrow \mathcal{G}$ between $\mathcal{E}$ and $\mathcal{G}$ is a morphism of the underlying sheaves $\varphi\colon  E \rightarrow G$ for which there is an element $\lambda\in k$ such that $\beta \circ \varphi=\lambda\alpha.$ We say that $\varphi\colon  \mathcal{E}\rightarrow \mathcal{G}$ is injective (resp.\ surjective) if the morphism $\varphi\colon  E \rightarrow G$ is injective (resp.\ surjective).
\end{definition}
\begin{remark}\label{rem:induced} 
Let $\mathcal{E}=(E, \alpha)$ be a framed sheaf. If $E'$ is a subsheaf of $E$ with quotient $E''=E/E'$, then we have the following commutative diagram
\begin{equation*}
  \begin{tikzpicture}
    \def\x{1.5}
    \def\y{-1.2}
    \node (A0_0) at (0*\x, 0*\y) {$0$};
    \node (A0_1) at (1*\x, 0*\y) {$E'$};
    \node (A0_2) at (2*\x, 0*\y) {$E$};
    \node (A0_3) at (3*\x, 0*\y) {$E''$};
    \node (A0_4) at (4*\x, 0*\y) {$0$};
    \node (A1_1) at (1*\x, 1*\y) {$F$};
    \node (A1_2) at (2*\x, 1*\y) {$F$};
    \node (A1_3) at (3*\x, 1*\y) {$F$};
    \path (A0_1) edge [->] node [auto] {$\scriptstyle{i}$} (A0_2);
    \path (A0_0) edge [->] node [auto] {$\scriptstyle{}$} (A0_1);
    \path (A0_3) edge [->] node [auto] {$\scriptstyle{\alpha''}$} (A1_3);
    \path (A0_2) edge [->] node [auto] {$\scriptstyle{\alpha}$} (A1_2);
    \path (A0_3) edge [->] node [auto] {$\scriptstyle{}$} (A0_4);
    \path (A1_1) edge [->] node [auto] {$\scriptstyle{\cdot \lambda}$} (A1_2);
    \path (A1_2) edge [->] node [auto] {$\scriptstyle{\cdot \mu}$} (A1_3);
    \path (A0_2) edge [->] node [auto] {$\scriptstyle{q}$} (A0_3);
    \path (A0_1) edge [->] node [auto] {$\scriptstyle{\alpha'}$} (A1_1);
  \end{tikzpicture}
\end{equation*}
where $\lambda=0, \mu=1$ if $\alpha'=0$, and $\lambda=1, \mu=0$ if $\alpha'\neq 0.$ Thus the inclusion morphism $i$ (resp.\ the projection morphism $q$) induces a morphism of framed sheaves between $\mathcal{E}'$ and $\mathcal{E}$ (resp.\ $\mathcal{E}$ and $\mathcal{E}/E'$). Note that in general an injective (resp.\ surjective) morphism $E\rightarrow G$ between the underlying sheaves of two framed sheaves $\mathcal{E}=(E,\alpha)$ and $\mathcal{G}=(G,\beta)$ does not lift to a morphism $\mathcal{E}\rightarrow \mathcal{G}$ of the corresponding framed sheaves. \triend
\end{remark}
\begin{lemma}{\normalfont (\cite[Lemma 1.5]{art:huybrechtslehn1995-II}).}
Let $\mathcal{E}=(E, \alpha)$ and $\mathcal{G}=(G,\beta)$ be framed sheaves. The set $\mathrm{Hom}(\mathcal{E}, \mathcal{G})$ of morphisms of framed sheaves is a linear subspace of $\mathrm{Hom} (E,G)$. If $\varphi\colon  \mathcal{E}\rightarrow \mathcal{G}$ is an isomorphism, then the factor $\lambda$ in the definition can be taken in $k^{*}$. In particular, the isomorphism $\varphi_0=\lambda^{-1}\varphi$ satisfies $\beta\circ \varphi_0=\alpha.$ Moreover, if $\mathcal{E}$ and $\mathcal{G}$ are isomorphic, then their framed Hilbert polynomials coincide.
\end{lemma}
\begin{proposition}\label{prop:quozim}
Let $\mathcal{E}=(E, \alpha)$ and $\mathcal{G}=(G,\beta)$ be framed sheaves. If $\varphi$ is a nonzero morphism of framed sheaves between $\mathcal{E}$ and $\mathcal{G}$, then
\begin{equation*}
P(\mathcal{E}/\ker\varphi)\leq P(\mathrm{Im}\:\varphi, \beta').
\end{equation*}
\end{proposition}
\proof
Consider a morphism of framed sheaves $\varphi\in\mathrm{Hom}(\mathcal{E},\mathcal{G})$, $\varphi\neq 0.$ There exists $\lambda \in k$ such that $\beta \circ \varphi=\lambda\alpha.$ Note that $E/\ker\varphi\simeq \mathrm{Im}\:\varphi$ hence their Hilbert polynomials coincide. It remains to prove that $\epsilon(\alpha'')\geq \epsilon(\beta').$ If $\lambda=0$, then $\beta'=0$ and therefore $\epsilon(\beta')=0\leq \epsilon(\alpha'').$ Assume now $\lambda\neq 0$: $\alpha=0$ if and only if $\beta\vert_{\mathrm{Im}\:\varphi}=0$, hence $\epsilon(\beta')=0=\epsilon(\alpha'').$ If $\alpha\neq 0$, then also $\alpha''\neq 0$. Indeed if $\alpha''=0$, then $\alpha\vert_{\ker\varphi}\neq 0$; this implies that $\lambda (\alpha\vert_{\ker\varphi})=(\beta \circ \varphi)\vert_{\ker\varphi}=0$ and therefore $\lambda=0$, but this is in contradiction with our previous assumption. Thus, if $\lambda\neq 0$ and $\alpha\neq 0$ then we obtain $\epsilon(\beta')=1=\epsilon(\alpha'').$
\qedhere\endproof
\begin{remark}
Let $\mathcal{E}=(E, \alpha)$ and $\mathcal{G}=(G,\beta)$ be framed sheaves and $\varphi\colon \mathcal{E}\rightarrow \mathcal{G}$ a nonzero morphism of framed sheaves. By the previous proposition, we get
\begin{equation*}
P(\mathcal{E})=P(\ker\varphi,\alpha')+P(\mathcal{E}/\ker\varphi)\leq P(\ker\varphi,\alpha')+P(\mathrm{Im}\:\varphi, \beta'). 
\end{equation*}
The inequality may be strict. This phenomenon does not appear in the nonframed case and depends on the fact that in general the isomorphism  $E/\ker\varphi\cong \mathrm{Im}\:\varphi$ does not induce an isomorphism $\mathcal{E}/\ker\varphi\cong (\mathrm{Im}\:\varphi, \beta')$ (there are examples where indeed it does not). \triend
\end{remark}

\section{Semistability}\label{sec:semistability}

In this section we give a generalization to framed sheaves of the Gieseker (semi)stability condition for $d$-dimensional sheaves (see \cite[Definition 1.2.4]{book:huybrechtslehn2010}). Comparing with the classical case, the (semi)stability condition for framed sheaves has an additional parameter $\delta$, which is a polynomial with rational coefficients. This definition was given in Huybrechts and Lehn's article \cite{art:huybrechtslehn1995-I} for a nonsingular polarized variety of dimension $d$; we generalize it to polarized schemes. We also need to apply this definition to sheaves of dimension smaller than $d$, due to the fact that even if we want to work only with framed sheaves $\mathcal{E}=(E,\alpha)$ with $E$ pure sheaf of dimension $d$, the graded factors of the \emph{framed} Harder-Narasimhan or Jordan-H\"older filtrations of $\mathcal{E}$ may have dimension less than $d.$ We will also present examples where the underlying sheaf of a semistable framed sheaf is not pure, and examples of non-semistable framed sheaves $(E,\alpha)$ with $E$ Gieseker semistable (see Example \ref{ex:p2}).

Recall that there is a natural ordering of rational polynomials given by the lexicographic order of their coefficients. Explicitly, $f\leq g$ if and only if $f(m)\leq g(m)$ for $m\gg 0$. Analogously, $f<g$ if and only if $f(m)<g(m)$ for $m\gg 0.$

We shall use the following convention: if the word ``(semi)\-stable'' occurs in any statement in combination with the symbol $(\leq)$, then two variants of the statement are asserted at the same time: a ``semistable'' one involving the relation ``$\leq$'' and a ``stable'' one involving the relation ``$<$''.

We now give a definition of (semi)stability for $d$-dimensional framed sheaves.
\begin{definition}\label{def:semi}
A $d$-dimensional framed sheaf $\mathcal{E}=(E, \alpha)$ is said to be \textit{(semi)stable} with respect to $\delta$ if and only if the following conditions are satisfied:
\begin{itemize}
\item[(a)] $\beta_d(E)P(E')\; (\leq)\; \beta_d(E') P(\mathcal{E})$ for all subsheaves $E'\subseteq \ker\alpha$,
\item[(b)] $\beta_d(E)(P(E')-\delta)\; (\leq)\; \beta_d(E') P(\mathcal{E})$ for all subsheaves $E'\subset E.$
\end{itemize}
\end{definition}
\begin{lemma}{\normalfont (\cite[Lemma 1.2]{art:huybrechtslehn1995-II}).}\label{lem:1.2}
Let $\mathcal{E}=(E,\alpha)$ be a $d$-dimensional framed sheaf. If $\mathcal{E}$ is semistable with respect to $\delta$, then $\ker\alpha$ is a pure sheaf of dimension $d.$
\end{lemma}
\proof
Let $T:=T_{d-1}(\ker\alpha).$ By the semistability condition (a), we get
\begin{equation*}
 \beta_d(E)P(T)\leq \beta_d(T) P(\mathcal{E}).
\end{equation*}
Since $\beta_d(T)=0$, we get $P(T)\; \leq\; 0.$ On the other hand, if $T\neq 0$, the leading coefficient of $P(T)$ is positive. Thus we get a contradiction and therefore $T=0.$
\qedhere\endproof
In the following, we shall call $d$\textit{-torsion sheaf} a sheaf $E$ of dimension less than $d.$ In this case, $T_{d-1}(E)=E.$
\begin{example}\label{ex:p2}
Let $(X,\mathcal{O}_X(1))$ be a nonsingular polarized variety of dimension $d$ and $D=D_1+\cdots +D_l$ an effective divisor on $X$, where $D_1, \ldots, D_l$ are distinct prime divisors. Consider the short exact sequence associated to the line bundle $\mathcal{O}_{X}(-D)$:
\begin{equation*}
 0\longrightarrow \mathcal{O}_{X}(-D)\longrightarrow \mathcal{O}_{X} \stackrel{\alpha}{\longrightarrow} i_{*}(\mathcal{O}_{Y})\longrightarrow 0,
\end{equation*}
where $Y:=\mathrm{Supp}(D)=D_1\cup \cdots \cup D_l.$ Note that $P(i_{*}(\mathcal{O}_{Y}))$ is a rational polynomial of degree $d-1.$ Let $\delta(n)\in \mathbb{Q}[n]$ be a polynomial of degree $d-1$ such that $\delta>P(i_{*}(\mathcal{O}_{Y})).$ Then we get
\begin{equation*}
 P(\mathcal{O}_{X})-\delta<P(\mathcal{O}_X)-P(i_{*}(\mathcal{O}_{Y}))=P(\mathcal{O}_{X}(-D))<P(\mathcal{O}_{X}).
\end{equation*}
Thus in this way we have that the $d$-dimensional framed sheaf $(\mathcal{O}_{X}, \alpha\colon \mathcal{O}_{X} \rightarrow i_*(\mathcal{O}_{Y}))$ is not semistable with respect to $\delta.$ We thus have obtained an example of a framed sheaf which is not semistable with respect to a fixed $\delta$ but the underlying sheaf is Gieseker semistable. It is possible to construct examples of semistable framed sheaves whose underlying sheaves are not Gieseker semistable, how we will see in Example \ref{ex:seminot}.

On the other hand, if we define $\alpha$ as the projection from $\mathcal{O}_X(-D)\oplus i_*(\mathcal{O}_Y)$ to its second factor, it is easy to check that the $d$-dimensional framed sheaf $(\mathcal{O}_X(-D)\oplus i_*(\mathcal{O}_Y), \alpha)$ is semistable with respect to $\delta:=P(i_*(\mathcal{O}_Y))$ and the underlying sheaf contains a nonzero $d$-torsion subsheaf.
\triend\end{example}
\begin{definition}
A $d$-dimensional framed sheaf $\mathcal{E}=(E,\alpha)$ is \textit{geometrically stable} with respect to $\delta$ if for any base extension $X\times_{\mathrm{Spec}(k)}\mathrm{Spec}(K)\stackrel{f}{\rightarrow} X$, the pull-back $f^*(\mathcal{E}):=(f^*(E),f^*(\alpha))$ is stable with respect to $\delta.$
\end{definition}
Since $X$ is defined over an algebraically closed field, in the unframed case Gieseker stability implies geometrical Gieseker stability (see \cite[Definition 1.2.9 and Corollary 1.5.11]{book:huybrechtslehn2010}). At the moment we do not know if the two definitions are equivalent in the framed case; we will only prove that they coincide for a particular class of $d$-dimensional framed sheaves (see Corollary \ref{cor:geostable}).

We have the following characterization of the semistability condition in terms of quotients:
\begin{proposition}\label{prop:quoziente}
Let $\mathcal{E}=(E,\alpha)$ be a $d$-dimensional framed sheaf. Then the following conditions are equivalent:
\begin{itemize}
\item[(a)] $\mathcal{E}$ is (semi)stable with respect to $\delta.$
\item[(b)] For any nonzero surjective morphism of framed sheaves $\varphi\colon \mathcal{E}\rightarrow (Q,\beta)$, one has\newline $\beta_d(Q)p(\mathcal{E})\;(\leq)\; P(Q,\beta).$
\end{itemize}
\end{proposition}
\proof
By using Proposition \ref{prop:quozim}, the assertion follows from the same arguments as in the nonframed case (see \cite[Proposition 1.2.6]{book:huybrechtslehn2010}).
\qedhere\endproof
In the papers by Huybrechts and Lehn, one finds two different definitions of (semi)stability of rank zero framed sheaves on nonsingular polarized varieties. In \cite{art:huybrechtslehn1995-I}, they use the same definition for the framed sheaves of positive or zero rank, and with that definition, all framed sheaves of rank zero are automatically semistable but not stable (with respect to any stability polynomial $\delta$). According to \cite[Definition 1.1]{art:huybrechtslehn1995-II}, the semistability of a rank zero framed sheaf depends on the choice of a stability polynomial $\delta$, but all semistable framed sheaves of rank zero are automatically stable. We give a new definition of the (semi)stability for $d$-torsion framed sheaves on polarized schemes which singles out exactly those objects which may appear as $d$-torsion components of the Harder-Narasimhan and Jordan-H\"older filtrations.
\begin{definition}\label{def:semi-rankzero}
Let $\mathcal{E}=(E,\alpha)$ be a $d$-torsion framed sheaf. If $\alpha$ is injective, we say that $\mathcal{E}$ is \textit{semistable}\footnote{For $d$-torsion sheaves, the definition of semistability of the corresponding framed sheaves does not depend on $\delta.$}. Moreover, if $P(E)=\delta$ we say that $\mathcal{E}$ is \textit{stable} with respect to $\delta.$
\end{definition}
We conclude this section by giving some constraints on the choice of the stability polynomial.
\begin{lemma}{\normalfont (\cite[Lemma 2.1]{art:huybrechtslehn1995-I}).}
Let $F$ be a sheaf of dimension less or equal to $d.$ Let $\mathcal{E}=(E,\alpha\colon E\rightarrow F)$ be a $d$-dimensional framed sheaf with $\ker\alpha$ nonzero and $\alpha$ surjective. If $\mathcal{E}$ is (semi)stable with respect to $\delta$, then
\begin{equation*}
\delta\; (\leq)\;  P(E)-\frac{\beta_d(E)}{\beta_d(\ker\alpha)}(P(E)-P(F)).
\end{equation*}
If $F$ is a $d$-torsion sheaf, $\delta\; (\leq)\;  P(F).$
\end{lemma}
\proof
By the (semi)stability condition (a), we get
\begin{equation*}
 \beta_d(E)P(\ker\alpha)\;(\leq)\; \beta_d(\ker\alpha)P(\mathcal{E})=\beta_d(\ker\alpha)\left(P(E)-\delta\right).
\end{equation*}
By Lemma \ref{lem:1.2} $\beta_d(\ker\alpha)>0$ , hence $\delta\; (\leq)\;  P(E)-(\beta_d(E)/\beta_d(\ker\alpha))P(\ker\alpha).$ Since $P(E)-P(\ker\alpha)=P(\mathrm{Im}\:\alpha)=P(F)$, we obtain the assertion. If $F$ is a $d$-torsion sheaf, $\beta_d(\mathrm{Im}\:\alpha)=0.$ Therefore $\beta_d(\ker\alpha)=\beta_d(E)$ and $\delta\;(\leq)\;P(E)-P(\ker\alpha)=P(F).$
\qedhere\endproof

\section{Characterization of semistability by means of framed saturated subsheaves}

Let $\mathcal{E}=(E,\alpha)$ be a $d$-dimensional framed sheaf, and assume that $\ker\alpha$ is a pure sheaf of dimension $d.$ In this section we would like to answer the following question: to verify if $\mathcal{E}$ is (semi)stable or not, do we need to check the inequalities (a) and (b) in Definition \ref{def:semi} for all subsheaves of $E$? Or, can we restrict our attention to a smaller family of subsheaves of $E$? For Gieseker (semi)stability condition, this latter family consists of \emph{saturated} subsheaves of $E$ (see \cite[Proposition 1.2.6]{book:huybrechtslehn2010}). In the framed case, we need to enlarge this family because of the framing, as we explain in what follows.
\begin{definition}
Let $E$ be a sheaf. The \textit{saturation} of a subsheaf $E'\subset E$ is the minimal subsheaf $\bar{E'}\subset E$ containing $E'$ such that the quotient $E/\bar{E'}$ is pure of dimension $\dim(E)$ or zero.
\end{definition}
Now we generalize this definition to framed sheaves:
\begin{definition}
Let $\mathcal{E}=(E,\alpha)$ be a $d$-dimensional framed sheaf where $\ker\alpha$ is a pure sheaf of dimension $d.$ Let $E'$ be a subsheaf of $E.$ The \textit{framed saturation} $\bar{E'}$ of $E'$ is the saturation of $E'$ as subsheaf of
\begin{itemize}
\item $\ker\alpha$, if $E'\subseteq \ker\alpha.$
\item $E$, if $E'\not\subseteq \ker\alpha.$
\end{itemize}
\end{definition}
\begin{remark}\label{rem:satura}
Let $\bar{E'}$ be the framed saturation of $E'\subset E.$ In the first case described in the definition, if $\beta_d(E')<\beta_d(\ker\alpha)$, the quotient $Q=E/\bar{E'}$ is a $d$-dimensional sheaf, with nonzero induced framing $\beta$, and fits into an exact sequence
\begin{equation}\label{eq:betano}
0\longrightarrow Q'\longrightarrow Q\stackrel{\beta}{\longrightarrow} \mathrm{Im}\:\alpha\longrightarrow 0,
\end{equation}
where $Q'=\ker\beta$ is a pure quotient of $\ker\alpha$ of dimension $d.$ If $\beta_d(E')=\beta_d(\ker\alpha)$, then $\bar{E'}=\ker\alpha$ and $Q=E/\ker\alpha.$ In the second case, if $\beta_d(E')=\beta_d(E)$, then $\bar{E'}=E$ and $Q=0.$ Otherwise, $Q$ is a pure sheaf of dimension $d$ with zero induced framing. Moreover $\beta_d(E')=\beta_d(\bar{E'})$, $P(E')\leq P(\bar{E'})$ and $P(\mathcal{E}')\leq P(\bar{\mathcal{E}'}).$ \triend
\end{remark}
\begin{example}
Let us consider the framed sheaf $(\mathcal{O}_{X}, \alpha\colon \mathcal{O}_{X} \rightarrow i_*(\mathcal{O}_{Y}))$ on $X$, defined in Example \ref{ex:p2}. Since $\ker\alpha=\mathcal{O}_{X}(-D)$, the saturation of $\mathcal{O}_{X}(-D)$ (as subsheaf of $\mathcal{O}_{X}$) is $\mathcal{O}_{X}$ but the framed saturation of $\mathcal{O}_{X}(-D)$ is $\mathcal{O}_{X}(-D).$
\triend\end{example}
We have the following characterization:
\begin{proposition}\label{prop:quoziente1}
Let $\mathcal{E}=(E,\alpha)$ be a $d$-dimensional framed sheaf where $\ker\alpha$ is a pure sheaf of dimension $d.$ Then the following conditions are equivalent:
\begin{itemize}
\item[(a)] $\mathcal{E}$ is (semi)stable with respect to $\delta.$
\item[(b)] For any framed saturated subsheaf $E'\subset E$ one has $P(E',\alpha')\;(\leq)\; \beta_d(E')p(\mathcal{E}).$
\item[(c)] For any nonzero surjective morphism of framed sheaves $\varphi\colon  \mathcal{E}\rightarrow (Q,\beta)$, where $\alpha=\beta\circ \varphi$ and $Q$ is one of the following:
\begin{itemize}
\item a $d$-dimensional sheaf with nonzero framing $\beta$ such that $\ker\beta$ is a pure sheaf of dimension $d$,
\item a pure sheaf of dimension $d$ with zero framing $\beta$,
\item $Q=E/\ker\alpha,$
\end{itemize}
one has $\beta_d(Q)p(\mathcal{E})\leq P(Q,\beta).$
\end{itemize}
\end{proposition}
\proof
The implication $(a)\Rightarrow (b)$ is obvious. By Remark \ref{rem:satura}, $P(\mathcal{E}')\leq P(\bar{\mathcal{E}'})\leq\beta_d(\bar{E'})p(\mathcal{E})=\beta_d(E')p(\mathcal{E})$, where $\bar{E'}$ is the framed saturation of $E'$, thus  $(b)\Rightarrow (a).$ Finally, the framed sheaf $\mathcal{Q}$ has the properties asserted in condition (c) if and only if $\ker\varphi$ is a framed saturated subsheaf of $\mathcal{E}$, hence $(b)\Longleftrightarrow (c).$
\qedhere\endproof
\begin{corollary}\label{cor:semist}
Let $\mathcal{E}=(E, \alpha)$ and $\mathcal{G}=(G,\beta)$ be $d$-dimensional framed sheaves with the same reduced framed Hilbert polynomial $p.$ 
\begin{itemize}
 \item[(a)] If $\mathcal{E}$ is semistable and $\mathcal{G}$ is stable, then any nonzero morphism $\varphi\colon  \mathcal{E}\rightarrow \mathcal{G}$ is surjective. 
\item[(b)] If $\mathcal{E}$ is stable and $\mathcal{G}$ is semistable, then any nonzero morphism $\varphi\colon  \mathcal{E}\rightarrow \mathcal{G}$ is injective.
\item[(c)] If $\mathcal{E}$ and $\mathcal{G}$ are stable, then any nonzero morphism $\varphi\colon  \mathcal{E}\rightarrow \mathcal{G}$ is an isomorphism. Moreover, in this case $\mathrm{Hom}(\mathcal{E}, \mathcal{G})\simeq k.$ If in addition $\alpha\neq 0$, or equivalently, $\beta\neq 0$, there is a unique isomorphism $\varphi_0$ with $\beta \circ \varphi_0=\alpha.$
\end{itemize}
\end{corollary}
\begin{definition}
Let $\mathcal{E}$ be a framed sheaf. We say that $\mathcal{E}$ is \textit{simple} if $\mathrm{End}(\mathcal{E})\simeq k.$
\end{definition}
As in the unframed case, a stable $d$-dimensional framed sheaf is simple.

\section{Maximal destabilizing framed subsheaf}

Let $\mathcal{E}=(E, \alpha)$ be a $d$-dimensional framed sheaf where $\ker\alpha$ is a pure sheaf of dimension $d.$ If $\mathcal{E}$ is not semistable with respect to $\delta$, then there exist destabilizing subsheaves of $\mathcal{E}.$ In this section we would like to find the maximal one (with respect to the inclusion) and show that it has some interesting properties. Because of the framing, it is possible that this subsheaf is the $d$-torsion subsheaf of $E$ or is not saturated and we want to emphasize that these types of situation are not possible in the nonframed case (see \cite[Lemma 1.3.5]{book:huybrechtslehn2010}).
\begin{proposition}\label{prop:dest}
Let $\mathcal{E}=(E, \alpha)$ be a $d$-dimensional framed sheaf where $\ker\alpha$ is a pure sheaf of dimension $d.$ If $\mathcal{E}$ is not semistable with respect to $\delta$, then there is a framed saturated subsheaf $G\subset E$ such that for any subsheaf $E'\subseteq E$ one has
\begin{equation*}
\beta_d(E')P(\mathcal{G})\geq \beta_d(G)P(\mathcal{E}')
\end{equation*}
and in case of equality, one has $E'\subset G$.

Moreover, the framed sheaf $\mathcal{G}$ is uniquely determined and is semistable with respect to $\delta.$
\end{proposition}
\begin{definition}
We call $G$ \textit{the maximal destabilizing framed subsheaf of} $\mathcal{E}.$
\end{definition}
\proof[Proof of Proposition \ref{prop:dest}]
On the set of nontrivial subsheaves of $E$ we define the following order relation $\preceq$: let $G_1$ and $G_2$ be nontrivial subsheaves of $E$, $G_1\preceq G_2$ if and only if $G_1 \subseteq G_2$ and $\beta_d(G_2)P(\mathcal{G}_1)\leq \beta_d(G_1)P(\mathcal{G}_2).$ Since any ascending chain of subsheaves stabilizes, for any subsheaf $E'$, there is a subsheaf $G'$ such that $E'\subseteq G'\subseteq E$ and $G'$ is maximal with respect to $\preceq.$

First, assume that there exists a $d$-torsion subsheaf $E'\subset E$ such that $E'$ destabilizes $\mathcal{E}.$ This means $P(\mathcal{E}')>0$, that is, $P(E')>\delta.$  Let us consider the subsheaf $T_{d-1}(E)$ of $E.$ Then $P(T_{d-1}(E))\geq P(E')>\delta.$ Moreover, $E'\preceq T_{d-1}(E)$ and there are no $d$-dimensional subsheaves $G\subset E$ such that $T_{d-1}(E)\preceq G.$ Indeed, should that be the case, by the definition of $\preceq$ we would obtain $P(T_{d-1}(E))-\delta\leq 0$, in contradiction with the previous inequality. Thus we choose $G:=T_{d-1}(E).$ Since $\ker\alpha\vert_G=0$, $\mathcal{G}$ is semistable.

From now on we assume that for every $d$-torsion subsheaf $E'\subset E$ we have $P(\mathcal{E}')\leq 0.$ Let $G\subset E$ be a $\preceq$-maximal $d$-dimensional subsheaf with minimal multiplicity $\beta_d(G)$ among all $\preceq$-maximal subsheaves. Suppose there exists a subsheaf $H\subset E$ with $\beta_d(H)p(\mathcal{G})< P(\mathcal{H}).$ By hypothesis we have $\beta_d(H)>0.$ From $\preceq$-maximality of $G$ we get $G\nsubseteq H$ (in particular $H\neq E$). Now we want to show that we can assume $H\subset G$ by replacing $H$ with $G\cap H.$ 

If $H\nsubseteq G$, then the morphism $\varphi\colon H\rightarrow E\rightarrow E/G$ is nonzero. Moreover $\ker\varphi=G\cap H.$ The sheaf $I=\mathrm{Im}\:\varphi$ is of the form $J/G$ with $G\subsetneq J\subset E$ and $\beta_d(J)>0$. By the $\preceq$-maximality of $G$ we have $p(\mathcal{J})<p(\mathcal{G}),$ hence we obtain
\begin{equation*}
P(\mathcal{I})=P(\mathcal{J})-P(\mathcal{G})< \beta_d(J)p(\mathcal{G})-\beta_d(G)p(\mathcal{G})=\beta_d(I)p(\mathcal{G}),
\end{equation*}
and therefore
\begin{equation}\label{eq:im-min}
P(\mathcal{I})<\beta_d(I)p(\mathcal{G}).
\end{equation}
Now we want to prove the following:

\emph{Claim}: The sheaf $G\cap H$ is a $d$-dimensional subsheaf of $E.$
\proof 
Assume that $G\cap H=0.$ In this case, we get $H\cong I$; moreover this isomorphism lifts to an isomorphism $\mathcal{H}\cong \mathcal{I}$ of the corresponding framed sheaves and therefore $\varphi$ lifts to a morphism of framed sheaves $\varphi$ between $\mathcal{H}$ and $\mathcal{E}/G.$ From formula \eqref{eq:im-min} it follows
\begin{equation*}
p(\mathcal{G})< p(\mathcal{H})=p(\mathcal{I})<p(\mathcal{G}),
\end{equation*}
which is absurd. 

The sheaf $G\cap H$ is $d$-dimensional, indeed if we assume that $\beta_d(G\cap H)=0$, then we have $\beta_d(I)=\beta_d(H)$ and again by Proposition \ref{prop:quozim} and formula \eqref{eq:im-min} we get
\begin{eqnarray*}
P(G\cap H,\alpha')&=&P(\mathcal{H})-P(\mathcal{H}/(G\cap H))\\
&\geq&P(\mathcal{H})-P(\mathcal{I})>P(\mathcal{H})-\beta_d(H)P(\mathcal{G})>0
\end{eqnarray*}
hence $G\cap H$ is a $d$-torsion subsheaf of $E$ with $P(G\cap H,\alpha')>0$, but this is in contradiction with the hypothesis.
\black

By the following computation:
\begin{eqnarray*}
\beta_d(G\cap H)\left(p(G\cap H,\alpha')-p(\mathcal{H})\right)&=&\beta_d(H/(G\cap H))p(\mathcal{H})-P(\mathcal{H}/(G\cap H))\\
&>&\beta_d(I)p(\mathcal{H})-P(\mathcal{I})>\beta_d(I)\left(p(\mathcal{H})-p(\mathcal{G})\right)> 0
\end{eqnarray*}
we get $p(\mathcal{H})< p(G\cap H,\alpha')$, hence from now on we can consider a $d$-dimensional subsheaf $H\subset G$ such that $H$ is $\preceq$-maximal in $G$ and $p(\mathcal{G})<p(\mathcal{H}).$

Let $H'\subset E$ be a sheaf that contains $H$ and is $\preceq$-maximal in $E$. In particular, one has $p(\mathcal{G})<p(\mathcal{H})\leq p(\mathcal{H}').$ By $\preceq$-maximality of $H$ and $G$, we have $H'\nsubseteq G$. Then the morphism $\psi\colon H'\rightarrow E\rightarrow E/G$ is nonzero and $H\subset \ker\psi.$ As before $p(\mathcal{H}')< p(\ker\psi,\alpha').$ Thus we have $H\subset H'\cap G=\ker\psi$ and $p(\mathcal{H})< p(\ker\psi,\alpha')$, hence $H\preceq \ker\psi.$ This contradicts the $\preceq$-maximality of $H$ in $G$. Thus for all subsheaves $H\subseteq E$, we have $\beta_d(H)p(\mathcal{G})\geq P(\mathcal{H}).$

If there is a $d$-torsion subsheaf $H\subset E$ such that $P(\mathcal{H})=0$ and $H\nsubseteq G$, then by using the same argument as before, we get $P(H\cap G,\alpha')>0$, but this is in contradiction with the hypothesis. So there are no $d$-torsion subsheaves $H\subset E$ such that $P(\mathcal{H})=0$ and $H\nsubseteq G.$ If there is a $d$-dimensional subsheaf $H\subset E$ such that $p(\mathcal{G})= p(\mathcal{H})$, then $H\subset G$. In fact, if $H \nsubseteq G$ then we can replace $H$ by $G\cap H$ and using the same argument as before we obtain $p(\mathcal{G})= p(\mathcal{H})<p(H\cap G, \alpha')$ and this is absurd.
\qedhere\endproof

\subsection*{Minimal destabilizing framed quotient}

Let $\mathcal{E}=(E, \alpha)$ be a $d$-dimensional framed sheaf where $\ker\alpha$ is a pure sheaf of dimension $d.$ Suppose that $\mathcal{E}$ is not semistable with respect to $\delta.$
\begin{assumption}
If $\mathrm{Im}\:\alpha$ is a $d$-torsion sheaf, we further assume that $\ker\alpha$ is not the maximal destabilizing framed subsheaf. \triend
\end{assumption}
Let $T_1$ be the set consisting of nontrivial quotients $E\stackrel{q}{\rightarrow} Q\rightarrow 0$ such that
\begin{itemize}
\item $Q$ is a pure sheaf of dimension $d$,
\item the induced framing on $\ker q$ is nonzero,
\item $p(\mathcal{Q})<p(\mathcal{E}).$
\end{itemize}

Let $T_2$ be the set consisting of nontrivial quotients $E\stackrel{q}{\rightarrow} Q\rightarrow 0$ such that 
\begin{itemize}
\item $Q$ is a $d$-dimensional sheaf,
\item the induced framing on $\ker q$ is zero,
\item As in Remark \ref{rem:satura}, $Q$ fits into an exact sequence of the form \eqref{eq:betano},
\item $p(\mathcal{Q})<p(\mathcal{E}).$
\end{itemize}

By Proposition \ref{prop:quoziente1}, the set $T_1\cup T_2$ is nonempty. For $i=1,2$ define an order relation on $T_i$ as follows: if $Q_1, Q_2\in T_i$, we say that $Q_1 \sqsubset Q_2$ if and only if $p(\mathcal{Q}_1)< p(\mathcal{Q}_2)$ or $\beta_d(Q_1)< \beta_d(Q_2)$ in the case $p(\mathcal{Q}_1)= p(\mathcal{Q}_2).$ Let $Q_{-}^i$ be a $\sqsubset$-minimal element in $T_i$, for $i=1,2.$ Define
\begin{equation*}
Q_{-}:= \left\{ \begin{array}{ll}
         Q_{-}^1 & \mbox{if $p(\mathcal{Q}_{-}^1)<p(\mathcal{Q}_{-}^2)$ or if $p(\mathcal{Q}_{-}^2)=p(\mathcal{Q}_{-}^1)$ and $\beta_d(Q_{-}^1)\leq \beta_d(Q_{-}^2)$},
\vspace{.3cm}\\
        Q_{-}^2 & \mbox{if $p(\mathcal{Q}_{-}^2)<p(\mathcal{Q}_{-}^1)$ or if $p(\mathcal{Q}_{-}^2)=p(\mathcal{Q}_{-}^1)$ and $\beta_d(Q_{-}^2)< \beta_d(Q_{-}^1)$}.\end{array} \right.
\end{equation*}
We call $Q_{-}$ the \textit{minimal destabilizing framed quotient}. By easy computations one can prove the following:
\begin{proposition}\label{prop:minquotient}
The sheaf $G:=\ker(E\rightarrow Q_{-})$ is the maximal destabilizing framed subsheaf of $\mathcal{E}.$
\end{proposition}
By the uniqueness of the maximal destabilizing framed subsheaf, $Q_{-}$ is unique.

\section{Harder-Narasimhan filtration}

In this section we construct the Harder-Narasimhan filtration for a $d$-dimensional framed sheaf $\mathcal{E}=(E,\alpha)$ with nonzero framing. For the unframed case, we refer to \cite[Definition 1.3.2 and Theorem 1.3.4]{book:huybrechtslehn2010}. We adapt the techniques used by Harder and Narasimhan in the case of vector bundles on curves (see \cite{art:hardernarasimhan1974}). When $\mathrm{Im}\:\alpha$ is a $d$-torsion sheaf, the multiplicity of $\ker\alpha$ is equal to the multiplicity of $E$ and because of this fact we get a more involved characterization of the Harder-Narasimhan filtration than in the nonframed case (as one can see in Proposition \ref{prop:condequi}). The characterization of the Harder-Narasimhan filtration when $\mathrm{Im}\:\alpha$ has dimension $d$ is similar to the nonframed case (see \cite[Theorem 1.3.4]{book:huybrechtslehn2010}).

Therefore we consider separately whether $\beta_d(\mathrm{Im}\:\alpha)$ is zero or positive. In the first case, we can have two types of $d$-torsion sheaves as graded factors of the Harder-Narasimhan filtration of $\mathcal{E}$: $T_{d-1}(E)\subset E$ and the quotient $E/\ker\alpha.$ In the second case, the only $d$-torsion sheaf that can appear as a graded factor of the Harder-Narasiham filtration is $T_{d-1}(E).$

Consider first the case $\beta_d(\mathrm{Im}\:\alpha)=0.$
\begin{definition}\label{def:hana}
Let $\mathcal{E}=(E, \alpha)$ be a $d$-dimensional framed sheaf where $\ker\alpha$ is a pure sheaf of dimension $d$ and $\mathrm{Im}\:\alpha$ is a $d$-torsion sheaf. A \textit{Harder-Narasimhan filtration} for $\mathcal{E}$ is an increasing filtration of framed saturated subsheaves
\begin{equation}\label{eq:hana}
\mathrm{HN}_{\bullet}(\mathcal{E}): 0=\mathrm{HN}_{0}(\mathcal{E})\subset \mathrm{HN}_{1}(\mathcal{E})\subset \cdots \subset \mathrm{HN}_{l}(\mathcal{E})=E
\end{equation}
which satisfies the following conditions
\begin{itemize}
	\item[(a)] the quotient sheaf $gr_i^{\mathrm{HN}}(\mathcal{E}):=\mathrm{HN}_{i}(\mathcal{E})/\mathrm{HN}_{i-1}(\mathcal{E})$ with the induced framing $\alpha_i$ is a semista\-ble framed sheaf with respect to $\delta$ for $i=1, 2, \ldots, l.$
	\item[(b)] The quotient $E/\mathrm{HN}_{i-1}(\mathcal{E})$ has dimension $d$, the kernel of the induced framing is a pure sheaf of dimension $d$ and the subsheaf $gr_i^{\mathrm{HN}}(\mathcal{E})$ is the maximal destabilizing framed subsheaf of $\mathcal{E}/\mathrm{HN}_{i-1}(\mathcal{E})$ for $i=1, 2, \ldots, l-1.$
\end{itemize}
\end{definition}
\begin{lemma}\label{lem:rzc}
Let $\mathcal{E}=(E, \alpha)$ be a $d$-dimensional sheaf where $\ker\alpha$ is a pure sheaf of dimension $d$ and $\mathrm{Im}\:\alpha$ is a $d$-torsion sheaf. Suppose that $\mathcal{E}$ is not semistable (with respect to $\delta$). Let $G$ be the maximal destabilizing framed subsheaf of $\mathcal{E}.$ If $G\neq \ker\alpha$, then for every $d$-torsion subsheaf $T$ of $E/G$, we get $P(\mathcal{T})\leq 0.$
\end{lemma}
\proof
If the quotient $E/G$ is a pure sheaf of dimension $d$, the condition is trivially satisfied. Otherwise let $T\subset E/G$ be a $d$-torsion subsheaf with $P(\mathcal{T})>0.$ The sheaf $T$ is of the form $E'/G$, where $G\subset E'$ and $\beta_d(E')=\beta_d(G)$, hence we obtain $p(\mathcal{E}')>p(\mathcal{G})$, therefore $E'$ contradicts the maximality of $G.$
\qedhere\endproof
\begin{theorem}\label{teo:hana}
Let $\mathcal{E}=(E, \alpha)$ is a $d$-dimensional sheaf where $\ker\alpha$ is a pure sheaf of dimension $d$ and $\mathrm{Im}\:\alpha$ is a $d$-torsion sheaf. Then there exists a unique Harder-Narasimhan filtration for $\mathcal{E}.$
\end{theorem}
\proof
\underline{Existence}. If $\mathcal{E}$ is a semistable framed sheaf with respect to $\delta$, we put $l=1$ and a Harder-Narasimhan filtration is
\begin{equation*}
\mathrm{HN}_{\bullet}(\mathcal{E}): 0=\mathrm{HN}_{0}(\mathcal{E})\subset \mathrm{HN}_{1}(\mathcal{E})=E
\end{equation*}
Else there exists a subsheaf $E_1\subset E$ such that $E_1$ is the maximal destabilizing framed subsheaf of $\mathcal{E}$. If $E_1=\ker\alpha$, a Harder-Narasimhan filtration is
\begin{equation*}
\mathrm{HN}_{\bullet}(\mathcal{E}): 0=\mathrm{HN}_{0}(\mathcal{E})\subset \ker\alpha\subset \mathrm{HN}_{2}(\mathcal{E})=E
\end{equation*}
Otherwise, by Lemma \ref{lem:rzc} $\left(E/E_1, \alpha''\right)$ is a $d$-dimensional framed sheaf with $\ker\alpha''$ a pure sheaf of dimension $d$ and no $d$-torsion destabilizing framed subsheaves. If $\mathcal{E}/E_1$ is a semistable framed sheaf, a Harder-Narasimhan filtration is
\begin{equation*}
\mathrm{HN}_{\bullet}(\mathcal{E}): 0=\mathrm{HN}_{0}(\mathcal{E})\subset E_1\subset \mathrm{HN}_{2}(\mathcal{E})=E
\end{equation*}
Else there exists a $d$-dimensional subsheaf $E_2'\subset E/E_1$ such that $E_2'$ is the maximal destabilizing framed subsheaf of $\mathcal{E}/E_1.$ We denote by $E_2$ its pre-image in $E$. Now we apply the previous argument to $E_2$ instead of $E_1.$ Thus we can iterate this procedure and we obtain a finite length increasing filtration of framed saturated subsheaves of $E$, which satisfies conditions (a) and (b).

\underline{Uniqueness}. The uniqueness of the Harder-Narasimhan filtration follows from the uniqueness of the maximal destabilizing framed subsheaf.
\qedhere\endproof
\begin{remark}
By construction, for $i>0$ at most one of the framings $\alpha_i$ is nonzero and all but possibly one of the factors $gr_i^{\mathrm{HN}}(\mathcal{E})$ are pure sheaves of dimension $d.$ In particular if $gr_1^{\mathrm{HN}}(\mathcal{E})$ is a $d$-torsion sheaf, $gr_1^{\mathrm{HN}}(\mathcal{E})=T_{d-1}(E)$ and $\alpha_1\neq 0$; if $gr_l^{\mathrm{HN}}(\mathcal{E})$ is a $d$-torsion sheaf, $gr_l^{\mathrm{HN}}(\mathcal{E})=E/\ker\alpha$ and $\alpha_l\neq 0.$
\triend\end{remark}
Now we want to relate condition (b) in Definition \ref{def:hana} with the framed Hilbert polynomials of the pieces of the Harder-Narasimhan filtration. In particular we get the following.
\begin{proposition}\label{prop:condequi}
Let $\mathcal{E}=(E, \alpha)$ be a $d$-dimensional sheaf where $\ker\alpha$ is a pure sheaf of dimension $d$ and $\mathrm{Im}\:\alpha$ is a $d$-torsion sheaf. Suppose there exists a filtration of the form \eqref{eq:hana} satisfying condition (a). Then condition (b) is equivalent to the following:
\begin{itemize}
	\item[(b')] the quotient $\left(E/\mathrm{HN}_{j}(\mathcal{E}), \alpha''\right)$ is a $d$-dimensional framed sheaf where $\ker\alpha''$ is a pure sheaf of dimension $d$ for $j=1,2, \ldots, l-2$, it has no $d$-torsion destabilizing framed subsheaves, and
\begin{equation}\label{eq:dise}
\beta_d(gr_{i+1}^{\mathrm{HN}}(\mathcal{E}))P(gr_i^{\mathrm{HN}}(\mathcal{E}), \alpha_i)>\beta_d(gr_{i}^{\mathrm{HN}}(\mathcal{E}))P(gr_{i+1}^{\mathrm{HN}}(\mathcal{E}), \alpha_{i+1})
\end{equation}
for $i=1, \ldots, l-1.$
\end{itemize}
\end{proposition}
\proof
We prove this result by means of arguments similar to those used in \cite[Lemma 1.3.8]{art:hardernarasimhan1974} to get the analogous result.
\qedhere\endproof
Now we turn to the case in which $\mathrm{Im}\:\alpha$ has dimension $d.$ First, we give the following definition.
\begin{definition}
Let $\mathcal{E}=(E, \alpha)$ be a $d$-dimensional sheaf where $\ker\alpha$ is a pure sheaf of dimension $d$ and $\mathrm{Im}\:\alpha$ has dimension $d.$ A \textit{Harder-Narasimhan filtration} for $\mathcal{E}$ is an increasing filtration of framed saturated subsheaves
\begin{equation*}
\mathrm{HN}_{\bullet}(\mathcal{E}): 0=\mathrm{HN}_{0}(\mathcal{E})\subset \mathrm{HN}_{1}(\mathcal{E})\subset \cdots \subset \mathrm{HN}_{l}(\mathcal{E})=E
\end{equation*}
which satisfies the following conditions
\begin{itemize}
\item[(a)] the quotient sheaf $gr_i^{\mathrm{HN}}(\mathcal{E}):=\mathrm{HN}_{i}(\mathcal{E})/\mathrm{HN}_{i-1}(\mathcal{E})$ with the induced framing $\alpha_i$ is a semistable framed sheaf with respect to $\delta$ for $i=1, 2, \ldots, l.$
\item[(b)] the quotient $\left(E/\mathrm{HN}_{i}(\mathcal{E}), \alpha''\right)$ is a $d$-dimensional framed sheaf where $\ker\alpha''$ is a pure sheaf of dimension $d$ for $i=1, \ldots, l-1$, it has no $d$-torsion destabilizing framed subsheaves, and
\begin{equation*}
\beta_d(gr_{i+1}^{\mathrm{HN}}(\mathcal{E}))P(gr_i^{\mathrm{HN}}(\mathcal{E}), \alpha_i)>\beta_d(gr_{i}^{\mathrm{HN}}(\mathcal{E}))P(gr_{i+1}^{\mathrm{HN}}(\mathcal{E}), \alpha_{i+1}).
\end{equation*}
\end{itemize}
\end{definition}
In this case one can prove results similar to those stated in Lemma \ref{lem:rzc}, Theorem \ref{teo:hana} and Proposition \ref{prop:condequi}. In particular we get the following:
\begin{theorem}
Let $\mathcal{E}=(E, \alpha)$ be a $d$-dimensional sheaf where $\ker\alpha$ is a pure sheaf of dimension $d$ and $\mathrm{Im}\:\alpha$ has dimension $d.$ Then there exists a unique Harder-Narasimhan filtration for $\mathcal{E}.$ 
\end{theorem}
\begin{remark}
Let $(X, \mathcal{O}_X(1))$ be a nonsingular polarized variety. In \cite{art:gomezsolszamora2011}, the authors proved that for a framed sheaf $(E, \alpha\colon E\rightarrow \mathcal{O}_X)$, with $E$ locally free, its Harder-Narasimhan filtration coincides with its \emph{Kempf filtration} coming from GIT theory. \triend
\end{remark}

We conclude this section by proving a result about the maximal destabilizing framed subsheaf of a framed sheaf. This holds for $\mathrm{Im}\:\alpha$ of any dimension.
\begin{definition}
Let $B$ be a pure sheaf on $X.$ The \textit{maximal reduced Hilbert polynomial} of $B$ is the reduced Hilbert polynomial of the maximal Gieseker destabilizing subsheaf of $B.$ We denote it by $p_{max}(B).$
\end{definition}
\begin{lemma}\label{lem:destquo}
Let $\mathcal{E}=(E,\alpha)$ be a semistable $d$-dimensional framed sheaf and $B$ a pure sheaf of dimension $d$ with zero framing. Suppose that $p(\mathcal{E}) > p_{max}(B).$ Then $\mathrm{Hom}(\mathcal{E},(B,0))=0.$
\end{lemma}
\proof
Let $\varphi\in \mathrm{Hom}(\mathcal{E},(B,0))$, $\varphi\neq 0.$ Let $j$ be minimal such that $\varphi(E)\subset \mathrm{HN}_j(B)$. Then there exists a nonzero morphism of framed sheaves $\bar{\varphi}\colon\mathcal{E}\rightarrow gr_j^{HN}(B).$ By Propositions \ref{prop:quozim} and \ref{prop:quoziente1} we get 
\begin{equation*}
p(\mathcal{E})\leq p(\mathcal{E}/\ker\bar{\varphi})\leq p(\mathrm{Im}\:\bar{\varphi})\leq p(gr_j^{HN}(B))\leq p_{max}(B)
\end{equation*}
and this is a contradiction with our assumption.
\qedhere\endproof
\begin{proposition}\label{prop:destquo}
Let $\mathcal{E}=(E, \alpha)$ be a $d$-dimensional sheaf where $\ker\alpha$ is a pure sheaf of dimension $d.$ Assume that $\mathcal{E}$ is not semistable with respect to $\delta.$ Let $G$ be the maximal destabilizing framed subsheaf of $\mathcal{E}.$ Then
\begin{equation*}
\mathrm{Hom}\left(\mathcal{G}, \mathcal{E}/G\right)=0.
\end{equation*}
\end{proposition}
\proof
We have to consider separately four different cases.

\emph{Case 1}: $G=\ker\alpha.$ In this case by definition of morphism of framed sheaves, we get $\mathrm{Hom}\left(\mathcal{G}, \mathcal{E}/G\right)=0.$

\emph{Case 2}: $\alpha\vert_G=0$ and $\beta_d(G)<\beta_d(\ker\alpha).$ In this case $\mathrm{Hom}\left(\mathcal{G}, \mathcal{E}/G\right)=\mathrm{Hom}(G, \ker\alpha/G).$ Recall that $G$ is a Gieseker semistable $d$-dimensional sheaf and $\ker\alpha/G$ is a pure sheaf of dimension $d$; moreover from the maximality of $G$ follows that $p_G>p(T/G)$ for all subsheaves $T/G\subset \ker\alpha/G$, hence $p(G)>p_{max}(\ker\alpha/G)$ and by \cite[Lemma 1.3.3]{book:huybrechtslehn2010} we obtain the assertion.

\emph{Case 3}: $\alpha\vert_G\neq 0$ and $\beta_d(G)>0.$ In this case $E/G$ is a pure sheaf of dimension $d$ and the induced framing is zero. From the maximality of $\mathcal{G}$ it follows that $p(\mathcal{G})>p(T/G)$ for all subsheaves $T/G\subset E/G$, so we can apply Lemma \ref{lem:destquo} and we get the assertion.

\emph{Case 4}: $G=T_{d-1}(E).$ Let $\varphi\colon T_{d-1}(E)\rightarrow E/T_{d-1}(E).$ Since $\beta_d(\mathrm{Im}\:\varphi)=0$ and $E/T_{d-1}(E)$ is a pure sheaf of dimension $d$, we have $\mathrm{Im}\:\varphi=0$ and therefore we obtain the assertion.
\qedhere\endproof

\section{Jordan-H\"older filtration}\label{sec:jh}

By analogy with the study of Gieseker semistable sheaves, we will define Jordan-H\"older filtrations for semistable $d$-dimensional framed sheaves. Because of the framing, one needs to use Lemma \ref{lem:fra} in the construction of the filtration. Moreover, in general we cannot extend the notions of \emph{socle} and the \emph{extended socle} for stable pure sheaves to the framed case, because, for example, the sum of two framed saturated subsheaves may not be framed saturated, hence we construct these objects only for a smaller family of framed sheaves having some extra properties.
\begin{definition}
Let $\mathcal{E}=(E,\alpha)$ be a semistable $d$-dimensional framed sheaf. A \textit{Jordan-H\"older filtration} of $\mathcal{E}$ is a filtration
\begin{equation*}
E_{\bullet}: 0=E_0\subset E_1\subset\cdots \subset E_l=E
\end{equation*}
such that all the factors $E_i/E_{i-1}$ together with the induced framings $\alpha_i$ are stable with framed Hilbert polynomial $P(E_i/E_{i-1},\alpha_i)=\beta_d(E_i/E_{i-1})p(\mathcal{E}).$
\end{definition}
\begin{proposition}{\normalfont (\cite[Proposition 1.13]{art:huybrechtslehn1995-II}).}
Jordan-H\"older filtrations always exist. The framed sheaf
\begin{equation*}
gr(\mathcal{E})=(gr(E),gr(\alpha))=\bigoplus_i (E_i/E_{i-1},\alpha_i)
\end{equation*}
does not depend on the choice of the Jordan-H\"older filtration.
\end{proposition}
\begin{remark}
By construction for $i>0$ all subsheaves $E_i$ are framed saturated and the framed sheaves $\mathcal{E}_i$ are semistable with framed Hilbert polynomial $\beta_d(E_i)p(\mathcal{E}).$ In particular $\mathcal{E}_1$ is a stable framed sheaf. Moreover at most one of the framings $\alpha_i$ is nonzero and all but possibly one of the factors $E_i/E_{i-1}$ are pure sheaves of dimension $d.$ \triend
\end{remark}
\begin{lemma}\label{lem:one}
Let $\mathcal{E}=(E,\alpha)$ be a semistable $d$-dimensional framed sheaf. Then there exists at most one subsheaf $E'\subset E$ such that $\alpha\vert_{E'}\neq 0$, $\mathcal{E}'$ is a stable framed sheaf and $P(\mathcal{E}')=\beta_d(E')p(\mathcal{E}).$
\end{lemma}
\proof
Suppose that there exist subsheaves $E_1$ and $E_2$ of $E$ such that $\alpha\vert_{E_i}\neq 0$, the framed sheaf $\mathcal{E}_i$ is stable (with respect to $\delta$) and $P(\mathcal{E}_i)=r_i p(\mathcal{E})$, where $r_i=\beta_d(E_i)$, for $i=1,2.$ So we have $P(E_i)=r_ip(\mathcal{E})+\delta$ for $i=1,2.$ Let $E_{12}=E_1\cap E_2.$ Suppose that $E_{12}\neq 0$ and $\alpha\vert_{E_{12}}\neq 0.$ Put $r_{12}=\beta_d(E_{12}).$ Since $\mathcal{E}_i$ is stable, $P(E_{12})-\delta<r_{12}p(\mathcal{E}).$ Consider the exact sequence
\begin{equation*}
0\longrightarrow E_{12}\longrightarrow E_1\oplus E_2 \longrightarrow E_1+E_2\longrightarrow 0.
\end{equation*}
The induced framing on $E_1+E_2$ by $\alpha$ is nonzero; we denote it by $\gamma.$
\begin{eqnarray*}
P(E_1+E_2)&=&P(E_1)+P(E_2)-P(E_{12})=r_1p(\mathcal{E})+\delta+r_2p(\mathcal{E})+\delta-P(E_{12})\\
&>&\beta_d(E_1+E_2)p(\mathcal{E})+\delta
\end{eqnarray*}
and therefore
\begin{equation*}
P(E_1+E_2,\gamma)=P(E_1+E_2)-\delta>\beta_d(E_1+E_2)p(\mathcal{E}),
\end{equation*}
but this is a contradiction, because $\mathcal{E}$ is semistable. Now consider the case $\alpha\vert_{E_{12}}=0.$ By similar computations, we obtain
\begin{equation*}
P(E_1+E_2,\gamma)=P(E_1+E_2)-\delta>\beta_d(E_1+E_2)p(\mathcal{E})+\beta_d(E_1+E_2)\delta>\beta_d(E_1+E_2)p(\mathcal{E}),
\end{equation*}
but this is absurd. Thus $E_{12}=0$ and therefore $E_1+E_2=E_1\oplus E_2.$ In this case we get
\begin{eqnarray*}
P(E_1+E_2,\gamma)&=&P(E_1+E_2)-\delta=P(E_1)+P(E_2)-\delta\\
&=&r_1 p(\mathcal{E})+\delta+r_2 p(\mathcal{E})+\delta-\delta\\
&=&\beta_d(E_1+E_2)p(\mathcal{E})+\delta>\beta_d(E_1+E_2)p(\mathcal{E}),
\end{eqnarray*}
but this is not possible.
\qedhere\endproof
\begin{remark}
Let $\mathcal{E}=(E,\alpha)$ be a semistable $d$-dimensional framed sheaf. If there exists $E'\subset E$ such that $\beta_d(E')=0$ and $P(E')=\delta$, then $E'=T_{d-1}(E).$ Indeed from $P(T_{d-1}(E))\geq P(E')$ follows that $P(T_{d-1}(E))\geq \delta.$ Since $\mathcal{E}$ is semistable, we have $P(T_{d-1}(E))=\delta$ and so $E'=T_{d-1}(E).$ \triend
\end{remark}
By using similar computations as before, one can prove:
\begin{lemma}\label{lem:somma}
Let $\mathcal{E}=(E,\alpha)$ be a semistable $d$-dimensional framed sheaf. Let $E_1$ and $E_2$ be two different subsheaves of $E$ such that 
$P(\mathcal{E}_i)=\beta_d(E_i)p(\mathcal{E})$ for $i=1,2.$ Then $P(E_1+E_2, \alpha')=\beta_d(E_1+E_2)p(\mathcal{E})$ and $P(E_1\cap E_2, \alpha')=
\beta_d(E_1\cap E_2)p(\mathcal{E}).$
\end{lemma}

\subsection{Framed sheaves that are locally free along the support of the framing sheaf}\label{sec:locfreesupport}

In this section we assume that $X$ is an integral scheme of dimension $d$ and $F$ is supported on a divisor $D$ and is a locally free $\mathcal{O}_{D}$-module.
\begin{definition}\label{def:locfree}
Let $\mathcal{E}=(E,\alpha)$ be a $d$-dimensional framed sheaf on $X.$ We say that $\mathcal{E}$ is $(D,F)$\textit{-framable} if $E$ is locally free in a neighbourhood of $D$ and $\alpha\vert_{D}$ is an isomorphism. We call $\mathcal{E}$ also $(D,F)$\textit{-framed sheaf}.
\end{definition}
Recall that in general for a $d$-dimensional framed sheaf $\mathcal{E}=(E,\alpha)$ where $\ker\alpha$ is torsion free, the subsheaf $T_{d-1}(E)$ of $E$ is supported on $\mathrm{Supp}(F).$ Therefore if $\mathcal{E}$ is $(D,F)$-framable, $E$ is torsion free.
\begin{example}\label{ex:seminot}
Let $\mathbb{CP}^2$ be the complex projective plane and $\mathcal{O}_{\mathbb{CP}^2}(1)$ the hyperplane line bundle. Let $l_{\infty}$ be a fixed line in $\mathbb{CP}^2$ and $i\colon l_{\infty}\rightarrow \mathbb{CP}^2$ the inclusion map. The torsion free sheaves of rank $r$ on $\mathbb{CP}^2$, trivial along a fixed line $l_{\infty}$ are --- in the language we introduced before --- $(l_{\infty}, \mathcal{O}_{l_{\infty}}^{\oplus r})$-framed sheaves of rank $r$ on $\mathbb{CP}^2.$ Let $\mathcal{M}(r,n)$ be the moduli space of $(l_{\infty}, \mathcal{O}_{l_{\infty}}^{\oplus r})$-framed sheaves of rank $r$ and second Chern class $n$ on $\mathbb{CP}^2.$ This moduli space is nonempty for $n\geq 1$ as one can see from the description of this space through ADHM data (see, e.g., \cite[Chapter 2]{book:nakajima1999}). Let $[(E,\alpha)]$ be a point in $\mathcal{M}(r,1)$: the sheaf $E$ is a torsion free sheaf of rank $r$ with second Chern class one. By \cite[Proposition 9.1.3]{book:lepotier1997}, $E$ is not Gieseker 
semistable. On the other hand, the framed sheaf $(E,\alpha)$ is stable with respect to a suitable choice of $\delta$ (see \cite[Theorem 3.1]{art:bruzzomarkushevich2011}). Thus we have proved that there exist semistable framed sheaves such that the underlying sheaves are not Gieseker semistable. \triend
\end{example}
\begin{lemma}
Let $\mathcal{E}=(E,\alpha)$ be a semistable $(D,F)$-framed sheaf. Let $E_1$ and $E_2$ be two different framed saturated subsheaves of $E$ such that 
$p(\mathcal{E}_i)=p(\mathcal{E})$, for $i=1,2.$ Assume that $\alpha\vert_{E_1}=0.$ Then $E_1+E_2$ is a framed saturated subsheaf of $E$ such that $gr(E_1+E_2,\alpha')=gr(\mathcal{E}_1)\oplus gr(\mathcal{E}_2).$ 
\end{lemma}
\proof
Since $\mathcal{E}$ is $(D,F)$-framable, the quotient $E/E_i$ is torsion free for $i=1,2$, hence $E/(E_1+E_2)$ is torsion free as well and therefore $E_1+E_2$ is framed saturated. 

By Lemma \ref{lem:somma}, $p(E_1+E_2,\alpha')=p(\mathcal{E}).$ Moreover we can always start with a Jordan-H\"older filtration of $\mathcal{E}_i$ and complete it to one of $(E_1+E_2,\alpha')$, hence we get $gr(\mathcal{E}_i)\subset gr(E_1+E_2,\alpha')$ (as framed sheaves) for $i=1,2.$ Let $G_\bullet: 0=G_0\subset G_1\subset\cdots \subset G_{l-1}\subset G_l=E_1$ be a Jordan-H\"older filtration for $\mathcal{E}_1$ and $H_\bullet: 0=H_0\subset H_1\subset \cdots \subset H_{s-1}\subset H_s=E_2$ a Jordan-H\"older filtration for $\mathcal{E}_2.$ Consider the filtration
\begin{equation*}
0=G_0\subset G_1\subset\cdots \subset G_{l-1}\subset G_l=E_1\subset E_1+H_p\subset \cdots \subset E_1+H_{t-1}\subset H_t=E_1+E_2
\end{equation*}
where $p=\min\{i\,\vert\, H_i\not\subset E_1\}.$ We want to prove that this is a Jordan-H\"older filtration for $(E_1+E_2,\alpha').$ It suffices to prove that $(E_1+H_j)/(E_1+H_{j-1})$ with its induced framing $\gamma_j$ is stable for $j=p,\ldots, t$ (we put $H_{p-1}=0$). First note that by Lemma \ref{lem:somma}, we get $P(E_1+H_j,\alpha')=\beta_d(E_1+H_j)p(\mathcal{E})$ and $P(E_1+H_{j-1},\alpha')=\beta_d(E_1+H_{j-1})p(\mathcal{E})$, hence
\begin{equation*}
P((E_1+H_j)/(E_1+H_{j-1}),\gamma_j)=\beta_d((E_1+H_j)/(E_1+H_{j-1}))p(\mathcal{E}).
\end{equation*}
Since $E/(E_1+H_{j-1})$ is torsion free, $\beta_d((E_1+H_j)/(E_1+H_{j-1}))>0.$ Let $T/(E_1+H_{j-1})$ be a subsheaf of $(E_1+H_{j})/(E_1+H_{j-1}).$ We have
\begin{eqnarray*}
&&P(T/(E_1+H_{j-1}),\gamma_j')=P(T,\alpha')-P(E_1+H_{j-1},\alpha')\leq \beta_d(T)p(\mathcal{E})-\beta_d(E_1+H_{j-1})p(\mathcal{E})\\
&&=\beta_d(T/(E_1+H_{j-1}))p(\mathcal{E})=\beta_d(T/(E_1+H_{j-1}))p((E_1+H_{j})/(E_1+H_{j-1}), \gamma_j),
\end{eqnarray*}
so the framed sheaf $((E_1+H_j)/(E_1+H_{j-1}),\gamma_j)$ is semistable. Moreover we can construct the following exact sequence of sheaves
\begin{equation*}
0\longrightarrow (E_1\cap H_j)/(E_1\cap H_{j-1})\longrightarrow H_j/H_{j-1}\stackrel{\varphi}{\longrightarrow} (E_1+H_j)/(E_1+H_{j-1})\longrightarrow 0.
\end{equation*}
Recall that the induced framing on $E_1$ is zero, hence the induced framing on $(E_1\cap H_j)/(E_1\cap H_{j-1})$ is zero and therefore the morphism $\varphi$ induces a surjective morphism between framed sheaves
\begin{equation*}
\varphi\colon  (H_j/H_{j-1},\beta_j)\longrightarrow ((E_1+H_j)/(E_1+H_{j-1}),\gamma_j).
\end{equation*}
Since $(H_j/H_{j-1},\beta_j)$ is stable, by Corollary \ref{cor:semist} the morphism $\varphi$ is injective, hence it is an isomorphism.
\qedhere\endproof
Now we introduce the \emph{extended framed socle} of a semistable $(D,F)$-framed sheaf, which plays a similar role of the maximal destabilizing framed subsheaf of a $d$-dimensional framed sheaf.
\begin{definition}
Let $\mathcal{E}=(E,\alpha)$ be a semistable $(D,F)$-framed sheaf. We call \textit{framed socle} of $\mathcal{E}$ the subsheaf of $E$ given by the sum of all framed saturated subsheaves $E'\subset E$ such that the framed sheaf $\mathcal{E}'=(E',\alpha\vert_{E'})$ is stable with reduced framed Hilbert polynomial $p(\mathcal{E}')=p(\mathcal{E}).$
\end{definition}
Let $\mathcal{E}=(E,\alpha)$ be a semistable $(D,F)$-framed sheaf. Consider the following two conditions on framed saturated subsheaves $E'\subset E$:
\begin{itemize}
\item[(a)] $p(\mathcal{E}')=p(\mathcal{E})$,
\item[(b)] each component of $gr(\mathcal{E}')$ is isomorphic (as a framed sheaf) to a subsheaf of $E.$
\end{itemize}
Let $E_1$ and $E_2$ two different framed saturated subsheaves of $E$ satisfying conditions (a) and (b). By previous lemmas the subsheaf $E_1+E_2$ is a framed saturated subsheaf of $E$ satisfying conditions (a) and (b) as well.
\begin{definition}
For a semistable $(D,F)$-framed sheaf $\mathcal{E}=(E,\alpha)$, we call \textit{extended framed socle} the maximal framed saturated subsheaf of $E$ satisfying the above conditions (a) and (b).
\end{definition}
\begin{proposition}\label{prop:extsocle}
Let $G$ be the extended framed socle of a semistable $(D,F)$-framed sheaf $\mathcal{E}=(E,\alpha).$ Then
\begin{itemize}
\item[(a)] $G$ contains the framed socle of $\mathcal{E}.$
\item[(b)] If $\mathcal{E}$ is simple and not stable, then $G$ is a proper subsheaf of $E.$
\end{itemize}
\end{proposition}
\proof
(a) It follows directly from the definition.

(b) Let $E_{\bullet}: 0=E_0\subset E_1\subset\cdots \subset E_l=E$ be a Jordan-H\"older filtration of $\mathcal{E}.$ If $E=G$, then the framed sheaf $(E/E_{l-1},\alpha_{l})$ is isomorphic (as framed sheaf) to a proper subsheaf $E'\subset E$ with induced framing $\alpha'.$ The composition of morphisms of framed sheaves
\begin{equation*}
  \begin{tikzpicture}
    \def\x{1.5}
    \def\y{-1.2}
    \node (A0_0) at (0*\x, 0*\y) {$E$};
    \node (A0_1) at (1*\x, 0*\y) {$E/E_{l-1}$};
    \node (A0_2) at (2*\x, 0*\y) {$E'$};
    \node (A0_3) at (3*\x, 0*\y) {$E$};
		\node (A1_0) at (0*\x, 1*\y) {$F$};
    \node (A1_1) at (1*\x, 1*\y) {$F$};
    \node (A1_2) at (2*\x, 1*\y) {$F$};
    \node (A1_3) at (3*\x, 1*\y) {$F$};
    \path (A0_0) edge [->] node [auto] {$\scriptstyle{\alpha}$} (A1_0);
		\path (A1_0) edge [->] node [auto] {$\scriptstyle{\cdot \nu}$} (A1_1);
		\path (A0_1) edge [->] node [auto] {$\scriptstyle{\sim}$} (A0_2);
    \path (A0_0) edge [->] node [auto] {$\scriptstyle{p}$} (A0_1);
    \path (A0_3) edge [->] node [auto] {$\scriptstyle{\alpha}$} (A1_3);
    \path (A0_2) edge [->] node [auto] {$\scriptstyle{\alpha'}$} (A1_2);
    \path (A1_1) edge [->] node [auto] {$\scriptstyle{\cdot \lambda}$} (A1_2);
    \path (A1_2) edge [->] node [auto] {$\scriptstyle{\cdot \mu}$} (A1_3);
    \path (A0_2) edge [->] node [auto] {$\scriptstyle{i}$} (A0_3);
    \path (A0_1) edge [->] node [auto] {$\scriptstyle{\alpha_l}$} (A1_1);
  \end{tikzpicture}
\end{equation*}
induces a morphism $\varphi\colon  \mathcal{E}\rightarrow \mathcal{E}$ that is not a scalar endomorphism of $\mathcal{E}.$ 
\qedhere\endproof
\begin{corollary}\label{cor:geostable}
A $(D,F)$-framed sheaf $\mathcal{E}=(E,\alpha)$ is stable with respect to $\delta$ if and only if it is geometrically stable. 
\end{corollary}
\proof
By using the previous proposition and the same arguments as in the unframed case (cf. \cite[Lemma 1.5.10 and Corollary 1.5.11]{book:huybrechtslehn2010}), we get the assertion.
\qedhere\endproof

\section{Relative Harder-Narasimhan filtration}\label{sec:relative}

In this section we construct a flat family of minimal destabilizing framed quotients associated to a flat family of $d$-dimensional framed sheaves. The construction in the framed case is somehow more complicated than in the nonframed case (see \cite[Theorem 2.3.2]{book:huybrechtslehn2010}), as one can see in what follows.

Let $g\colon Y\rightarrow S$ be a morphism of finite type of Noetherian schemes.
\begin{definition}
A \textit{flat family of sheaves on the fibres of the morphism} $g$ is a sheaf $G$ on $Y$, which is flat over $S.$
\end{definition}
Let $(f\colon \mathcal{X}\rightarrow S, \mathcal{O}_\mathcal{X}(1))$ be a relative polarized scheme of fixed dimension. Fix a positive integer $d$ such that $d\leq \dim(\mathcal{X}_s)$ for any $s\in S.$ Furthermore, fix a flat family of sheaves $F$ on the fibres of $f$ and a rational polynomial $\delta$ of degree $d-1$ and positive leading coefficient $\delta_1.$
\begin{definition}
A \textit{flat family of} $d$\textit{-dimensional framed sheaves on the fibres of the morphism} $f$ consists of a framed sheaf $\mathcal{E}=(E,\alpha\colon E\rightarrow F)$ on $\mathcal{X}$, where $E$ and $\mathrm{Im}\:\alpha$ are flat families of sheaves on the fibres of $f$, $\alpha_s\neq 0$ and $\dim(E_s)=d$ for all $s\in S.$
\end{definition}
\begin{remark}
By flatness of $E$ and $\mathrm{Im}\:\alpha$, we have that also $\ker\alpha$ is $S$-flat. \triend
\end{remark}
From now on, we fix $(f\colon \mathcal{X}\rightarrow S, \mathcal{O}_\mathcal{X}(1))$, $F$ and $\delta$ as introduced above, unless otherwise stated.

The direction we choose to obtain a flat family of minimal destabilizing framed quotients is the following: first we construct a \textit{universal quotient} (with fixed Hilbert polynomial) such that the induced framing is either nonzero at each fibre or zero at a generic fibre. In this way, generically, not only the Hilbert polynomial of that quotient is constant along the fibres, but also its framed Hilbert polynomial. Later we need to find a numerical polynomial such that the universal quotient with this polynomial as Hilbert polynomial gives the minimal destabilizing framed quotient at each fibre.

\subsection*{Relative framed Quot functor}

Let $(f\colon \mathcal{X}\rightarrow S, \mathcal{O}_\mathcal{X}(1))$, $F$ and $\delta$ be as introduced before.

Let $\mathcal{E}=(E, \alpha)$ be a flat family of $d$-dimensional framed sheaves on the fibres of $f$ and $P(n)\in \mathbb{Q}[n]$ a numerical polynomial. Define the contravariant functor from the category of Noetherian $S$-schemes of finite type to the category of sets
\begin{equation*}
\underline{\mathrm{FQuot}}_{\mathcal{X}/S}^P(\mathcal{E})\colon (Sch/S)\longrightarrow (Sets)
\end{equation*}
in the following way:
\begin{itemize}
\item For an object $T\rightarrow S$, $\underline{\mathrm{FQuot}}_{\mathcal{X}/S}^P(\mathcal{E})(T\rightarrow S)$ is the set consisting of the quotients (modulo isomorphism) $q\colon E_T\rightarrow Q$ such that
\begin{itemize}
\item[(a)] $Q$ is $T$-flat,
\item[(b)] the Hilbert polynomial of $Q_t$ is $P$ for all $t\in T$,
\item[(c)] there is an induced morphism $\widetilde{\alpha}\colon Q\rightarrow F_T$ such that $\widetilde{\alpha}\circ q=\alpha_T.$
\end{itemize}
\item For an $S$-morphism $g\colon  T'\rightarrow T$, $\underline{\mathrm{FQuot}}_{\mathcal{X}/S}^P(\mathcal{E})(g)$ is the map that sends $E_T\rightarrow Q$ to $E_{T'}\rightarrow g_{\mathcal{X}}^*Q$, where $g_{\mathcal{X}}\colon X_{T'}\rightarrow X_T$ is the induced morphism by $g.$
\end{itemize}

This functor is a subfunctor of the \textit{relative Quot functor} $\underline{\mathrm{Quot}}_{\mathcal{X}/S}^P(E)$, that is representable by a projective $S$-scheme $\pi\colon  \mathrm{Quot}_{\mathcal{X}/S}^P(E)\rightarrow S$ (cf. \cite[Theorem 2.2.4]{book:huybrechtslehn2010}).

The property (c) in the definition is closed, indeed by using the same arguments as in the proof of \cite[Theorem 1.6]{book:sernesi1986}, in particular results (iii) and (iv), one can see that property (c) is equivalent to the vanishing of some regular functions on $\mathrm{Quot}_{\mathcal{X}/S}^P(E).$ Hence, the functor $\underline{\mathrm{FQuot}}_{\mathcal{X}/S}^P(\mathcal{E})$ is representable by a closed subscheme $\mathrm{FQuot}_{\mathcal{X}/S}^P(\mathcal{E})\subset \mathrm{Quot}_{\mathcal{X}/S}^P(E).$ We denote by $\pi_{fr}$ the composition
\begin{equation*}
\pi_{fr}\colon \mathrm{FQuot}_{\mathcal{X}/S}^P(\mathcal{E})\hookrightarrow \mathrm{Quot}_{\mathcal{X}/S}^P(E) \stackrel{\pi}{\longrightarrow} S.
\end{equation*}
Roughly speaking, $\mathrm{FQuot}_{\mathcal{X}/S}^P(\mathcal{E})$ parametrizes all the quotients $E_s\stackrel{q}{\rightarrow} Q$, for $s\in S$, such that the induced framing on $\ker q$ is zero.

The \textit{universal object} on $\mathrm{FQuot}_{\mathcal{X}/S}^P(\mathcal{E})\times_S \mathcal{X}$ is the pull-back of the universal object on $\mathrm{Quot}_{\mathcal{X}/S}^P(E)\times_S \mathcal{X}$ with respect to the morphism $\mathrm{FQuot}_{\mathcal{X}/S}^P(\mathcal{E})\times_S \mathcal{X}\rightarrow \mathrm{Quot}_{\mathcal{X}/S}^P(E)\times_S \mathcal{X}$, induced by the closed embedding $\mathrm{FQuot}_{\mathcal{X}/S}^P(\mathcal{E})\hookrightarrow \mathrm{Quot}_{\mathcal{X}/S}^P(E).$

Let $s\in S$ and $q\in \pi_{fr}^{-1}(s)$ be $k$-rational points corresponding to the commutative diagram on $X_s$
\begin{equation*}
  \begin{tikzpicture}
    \def\x{1.5}
    \def\y{-1}
    \node (A0_0) at (0*\x, 0*\y) {$0$};
    \node (A0_1) at (1*\x, 0*\y) {$K$};
    \node (A0_2) at (2*\x, 0*\y) {$E_s$};
    \node (A0_3) at (3*\x, 0*\y) {$Q$};
    \node (A0_4) at (4*\x, 0*\y) {$0$};
    \node (A1_2) at (2*\x, 1*\y) {$F_s$};
  
    \path (A0_0) edge [->] node [auto] {$\scriptstyle{}$} (A0_1);
    \path (A0_1) edge [->] node [auto] {$\scriptstyle{i}$} (A0_2);
    \path (A0_3) edge [->] node [auto] {$\scriptstyle{\tilde{\alpha}}$} (A1_2);
    \path (A0_2) edge [->] node [left] {$\scriptstyle{\alpha_s}$} (A1_2);
    \path (A0_3) edge [->] node [auto] {$\scriptstyle{}$} (A0_4);
    \path (A0_2) edge [->] node [auto] {$\scriptstyle{q}$} (A0_3);
  \end{tikzpicture}
\end{equation*}
One has the following result about the tangent space of $\pi_{fr}^{-1}(s)$ at $q$:
\begin{proposition}\label{prop:tangentefibra}
The kernel of the linear map $(\mathrm{d}\pi_{fr})_q\colon T_q \mathrm{FQuot}_{\mathcal{X}/S}^P(\mathcal{E})\rightarrow T_s S$ is isomorphic to the linear space $\mathrm{Hom}(K, \ker\alpha_s/K)=\mathrm{Hom}(\mathcal{K},\mathcal{Q}).$
\end{proposition}
\proof
It suffices to use the same techniques of the proof of the corresponding result for $\pi$ (see \cite[Proposition 4.4.4]{book:sernesi2006}).
\qedhere\endproof
Now we have a tool for constructing a flat family of quotients (with a fixed Hilbert polynomial) of $\mathcal{E}$ such that the induced framing is nonzero in each fibre. Using the relative Quot scheme $\mathrm{Quot}_{\mathcal{X}/S}^P(E)$ associated to $E$, one can construct a flat family of quotients such that the induced framing is \textit{generically} zero.

\subsection*{Boundedness result}

Let $Y\rightarrow S$ be a projective morphism of Noetherian schemes and denote by $\mathcal{O}_Y(1)$ a line bundle on $Y$, which is very ample relative to $S.$

\begin{definition}
A family $\mathscr{G}$ of isomorphism classes of sheaves on the fibres of the morphism $Y\rightarrow S$ is \textit{bounded} if there is an $S$-scheme $T$ of finite type and a $\mathcal{O}_{Y_T}$-sheaf $E$ such that the given family is contained in the set $\{E_t\,\vert\, t \mbox{ is a closed point in } T\}.$
\end{definition}
Now we recall some characterizations of the property of boundedness which will be useful in the sequel.
\begin{lemma}{\normalfont (\cite[Lemma 2.5]{art:grothendieck1995-SB6}).}\label{lem:boun3}
Let $L$ be a sheaf on $Y$ and $\mathscr{G}$ the set of isomorphism classes of quotient sheaves $G$ of $L_s$ for $s$ running over the points of $S.$ Suppose that the dimension of $Y_s$ is $\leq r$ for all $s.$ Then the coefficient $\beta_r(G)$ is bounded from above and from below, and $\beta_{r-1}(G)$ is bounded from below. If $\beta_{r-1}(G)$ is bounded from above, then the family of sheaves $G/T_{\dim(G)-1}(G)$ is bounded.
\end{lemma}
\begin{corollary}\label{cor:boun-2}
Let $E$ be a flat family of sheaves on the fibres of a projective morphism $Y\rightarrow S.$ Then the family of pure quotients $Q$ of $E_s$ for $s\in S$ with hat-slopes bounded from above is a bounded family.
\end{corollary}
\begin{theorem}{\normalfont (\cite[Theorem 2.1]{art:grothendieck1995-SB6}).}\label{thm:bound-charact2}
The following properties of a family  $\mathscr{G}$ of isomorphism classes of sheaves on the fibres of $Y\rightarrow S$ are equivalent:
\begin{itemize}
 \item The family is bounded.
\item The set of Hilbert polynomials $\{P(G)\}_{G\in \mathscr{G}}$ is finite and there is a sheaf $E$ on $Y$ such that all $G\in \mathscr{G}$ admit surjective morphisms $E_s\rightarrow G.$
\end{itemize}
\end{theorem}
From this result and Corollary \ref{cor:boun-2} it follows that there are only a finite number of rational polynomials corresponding to Hilbert polynomials of destabilizing quotients $Q$ of $E_s$ for $s\in S.$ Thus it is possible to find the ``minimal'' polynomial that will be the Hilbert polynomial of the minimal destabilizing quotient of $E_s$ for a generic point $s\in S.$ Now we want to use the same argument in the framed case.

Let $(f\colon \mathcal{X}\rightarrow S, \mathcal{O}_\mathcal{X}(1))$, $F$ and $\delta$ be as before. Let $\mathcal{E}=(E, \alpha)$ be a flat family of $d$-dimensional framed sheaves on the fibres of $f.$

Let $\mathscr{F}_1$ be a family of quotients $E_s\stackrel{q}{\rightarrow} Q$, for $s\in S$, such that
\begin{itemize}
\item $\ker\alpha_s$ is a pure sheaf of dimension $d$,
\item $\ker q \not\subseteq \ker\alpha_s$,
\item $Q$ is a pure sheaf of dimension $d$ and $\hat{\mu}(Q)< \hat{\mu}(E_s).$
\end{itemize}
\begin{proposition}\label{prop:bounded1}
The family $\mathscr{F}_1$ is bounded.
\end{proposition}
\proof
The family $\mathscr{F}_1$ is contained in the family of purely $d$-dimensional quotients of $E$, with hat-slopes bounded from above, hence it is bounded by Corollary \ref{cor:boun-2} and Theorem \ref{thm:bound-charact2}.
\qedhere\endproof
Let $\mathscr{F}_2$ be a family of quotients
\begin{equation*}
  \begin{tikzpicture}
    \def\x{1.8}
    \def\y{-1.5}
    \node (A0_2) at (2*\x, 0*\y) {$E_s$};
    \node (A0_3) at (3*\x, 0*\y) {$Q$};
    \node (A0_4) at (4*\x, 0*\y) {$0$};
    \node (A1_2) at (2*\x, 1*\y) {$F_s$};
    \path (A0_2) edge [->] node [auto] {$\scriptstyle{q}$} (A0_3);
    \path (A0_3) edge [->] node [auto] {$\scriptstyle{}$} (A0_4);
    \path (A0_2) edge [->] node [auto] {$\scriptstyle{\alpha_s}$} (A1_2);
		\path (A0_3) edge [->] node [auto] {$\scriptstyle{\tilde{\alpha}}$} (A1_2);
  \end{tikzpicture}
\end{equation*}	
for which 
\begin{itemize}
\item $\ker\alpha_s$ is a pure sheaf of dimension $d$,
\item $Q$ fits into a exact sequence
\begin{equation*}
0\longrightarrow Q' \longrightarrow Q\stackrel{\tilde{\alpha}}{\longrightarrow} \mathrm{Im}\:\alpha_s\longrightarrow 0
\end{equation*}
where $Q'=\ker\tilde{\alpha}$ is a nonzero purely $d$-dimensional quotient of $\ker\alpha_s$,
\item $\hat{\mu}(Q)< \hat{\mu}(E_s)+\delta_1.$
\end{itemize}
\begin{proposition}\label{prop:bounded2}
The family $\mathscr{F}_2$ is bounded.
\end{proposition}
\proof
Since a family given by extensions of elements from two bounded families is bounded (cf. \cite[Proposition 1.2]{art:grothendieck1995-SB6}), it suffices to prove that every element in $\mathscr{F}_2$ is an extension of two elements that belong to two bounded families. By definition of flat family of framed sheaves, the families $\{\ker\alpha_s\}_{s\in S}$ and $\{\mathrm{Im}\:\alpha_s\}_{s\in S}$ are bounded. 

So it remains to prove that the family of quotients $Q'$ is bounded. Since the family $\{\ker\alpha_s\}$ is bounded, there exists a sheaf $G$ on $\mathcal{X}$ such that we have surjective morphisms $G_s\rightarrow \ker\alpha_s$ (see Theorem \ref{thm:bound-charact2}), hence the compositions $G_s\rightarrow \ker\alpha_s\rightarrow Q'$ are surjective as well.

By Lemma \ref{lem:boun3}, the coefficient $\beta_{d}(Q')$ is bounded from above and from below and the coefficient $\beta_{d-1}(Q')$ is bounded from below, hence $\hat{\mu}(Q')$ is bounded from below.

Moreover, since $\{E_s\}$ is a bounded family, the coefficients of their Hilbert polynomials are uniformly bounded from above and from below, hence $\hat{\mu}(E_s)$ is uniformly bounded from above and from below. Therefore, from the inequality $\hat{\mu}(Q)< \hat{\mu}(E_s)+\delta_1$, it follows that $\hat{\mu}(Q)$ is uniformly bounded from above. By using the fact that also $\{\mathrm{Im}\:\alpha_s\}$ is a bounded family, we obtain that $\hat{\mu}(Q')\leq A\hat{\mu}(Q)+B$ for some positive constants $A, B$, hence we get that $\hat{\mu}(Q')$ is uniformly bounded from above. 

Altogether, $\hat{\mu}(Q')$ is uniformly bounded from above and below, and by Lemma \ref{lem:boun3}, the family of quotients $Q'$ is bounded.
\qedhere\endproof

\subsection*{Relative minimal destabilizing framed quotient and Harder-Narasimhan filtration}

By using the same arguments as in the nonframed case (see \cite[Proposition 2.3.1]{book:huybrechtslehn2010}), we can prove the following:
\begin{proposition}\label{prop:open}
Let $(f\colon \mathcal{X}\rightarrow S, \mathcal{O}_\mathcal{X}(1))$, $F$ and $\delta$ be as in the assumptions of this section. Let $\mathcal{E}=(E, \alpha)$ be a flat family of $d$-dimensional framed sheaves on the fibres of $f.$ Assume that if $\beta_d(\mathrm{Im}\:\alpha_s)=0$, then $P(\mathrm{Im}\:\alpha_s)\geq \delta$ for any $s\in S.$ The set of points $s\in S$ such that $(E_s,\alpha_s)$ is (semi)stable with respect to $\delta$ is open in $S.$ 
\end{proposition}
Now we can prove the following generalization to the relative case of the theorem of the minimal destabilizing framed quotient.
\begin{theorem}\label{thm:rel}
Let $(f\colon \mathcal{X}\rightarrow S, \mathcal{O}_\mathcal{X}(1))$, $F$ and $\delta$ be as in the assumptions of this section. Let $\mathcal{E}=(E, \alpha)$ be a flat family of $d$-dimensional framed sheaves on the fibres of $f$ such that $\ker\alpha_s$ is a pure sheaf of dimension $d$ for a general point $s\in S.$ Assume that if $\beta_d(\mathrm{Im}\:\alpha_s)=0$, then $P(\mathrm{Im}\:\alpha_s)\geq \delta$ for any $s\in S.$ Then there is an integral $k$-scheme $T$ of finite type, a projective birational morphism $g\colon T\rightarrow S$, a dense open subscheme $U\subset T$ and a flat quotient $Q$ of $E_T$ such that for all points $t$ in $U$, $(E_t, \alpha_t)$ is a $d$-dimensional framed sheaf with $\ker\alpha_t$ pure sheaf of dimension $d$ and $\mathcal{Q}_t$ is the minimal destabilizing framed quotient of $\mathcal{E}_t$ with respect to $\delta$ or $Q_t=E_t.$ 

Moreover, the pair $(g, Q)$ is universal in the sense that if $\bar{g}\colon \bar{T}\rightarrow S$ is any dominant morphism of $k$-integral schemes and $\bar{Q}$ is a flat quotient of $E_{\bar{T}}$ satisfying the same property as $Q$, then there is a unique $S$-morphism $h\colon \bar{T}\rightarrow T$ such that $h_{\mathcal{X}}^{*}(Q)=\bar{Q}.$
\end{theorem}
\proof
Let $P$ be the Hilbert polynomial of $E$ and $r$ its leading coefficient. Denote by $p$ the rational polynomial $P/r.$

For $i=1,2$ let $A_i\subset \mathbb{Q}[n]$ be the set consisting of polynomials $P''$ such that there is a point $s\in S$ and a surjection $E_s\rightarrow E''$, where $P_{E''}=P''$ and $E''$ belongs to the family $\mathscr{F}_i$ (introduced in the previous section). By Propositions \ref{prop:bounded1} and \ref{prop:bounded2} and Theorem \ref{thm:bound-charact2}, the sets $A_1$ and $A_2$ are finite. Denote by $r''$ the leading coefficient of $P''$ and by $p''$ the rational polynomial $P''/r''$ (and similarly for other polynomials). Let
\begin{eqnarray*}
B_1&=&\left\{P''\in A_1\,\vert\, p''<p-\frac{\delta}{r}\right\},\\
B_2&=&\left\{P''\in A_2\,\vert\, p''-\frac{\delta}{r''}\leq p-\frac{\delta}{r}\right\}.
\end{eqnarray*}
The set $B_1\cup B_2$ is nonempty. Let $C_1$ be the set of polynomials $P''\in B_1$ such that $\pi(\mathrm{Quot}_{\mathcal{X}/S}^{P''}(E))=S$ and for any $s\in S$ one has $\pi^{-1}(s)\not\subset \mathrm{FQuot}_{\mathcal{X}/S}^{P''}(\mathcal{E}).$ Let $C_2$ be the set of polynomials $P''\in B_2$ such that $\pi_{fr}(\mathrm{FQuot}_{\mathcal{X}/S}^{P''}(\mathcal{E}))=S.$ Note that $C_1\cup C_2$ is nonempty. Now we want to find a polynomial $P_{-}$ in $C_1\cup C_2$ that is the Hilbert polynomial of the minimal destabilizing framed quotient of $E_s$ for a general point $s\in S.$

We define strict total order relations on $B_1$ and $B_2.$ In the first case, we set $P_1 \sqsubset P_2$ if and only if $p_1< p_2$ or $r_1< r_2$ in the case $p_1= p_2.$ In the second case, we stipulate that $P_1 \sqsubset P_2$ if and only if $p_1-\frac{\delta}{r_1}< p_2-\frac{\delta}{r_2}$ and $r_1< r_2$ in the case of equality. Let $P_{-}^i$ be the $\sqsubset$-minimal polynomial among all polynomials of $C_i$ for $i=1,2.$ Denote by $r^i_{-}$ the leading coefficient of $P^i_{-}$ and by $p_{-}^i$ the polynomial $P^i_{-}/r^i_{-}.$

Consider the following cases:
\begin{itemize}
\item \emph{Case 1}: $p^1_{-}<p^2_{-}-\frac{\delta}{r^2_{-}}.$ Put $P_{-}:=P_{-}^1.$ 
\vspace{0.2cm}
\item \emph{Case 2}: $p^1_{-}>p^2_{-}-\frac{\delta}{r^2_{-}}.$ Put $P_{-}:=P_{-}^2.$
\vspace{0.2cm}
\item \emph{Case 3}: $p^1_{-}=p^2_{-}-\frac{\delta}{r^2_{-}}.$ If $r_{-}^2< r_{-}^1$, put $P_{-}:=P_{-}^2$, otherwise $P_{-}:=P_{-}^1.$
\end{itemize}
Note that the set 
\begin{equation*}
\Big(\bigcup_{P''\in B_1\atop P''\sqsubset P_{-}^1} \pi(\mathrm{Quot}_{\mathcal{X}/S}^{P''}(E))\Big)\cup\Big(\bigcup_{P''\in B_2\atop P''\sqsubset P_{-}^2} \pi_{fr}(\mathrm{FQuot}_{\mathcal{X}/S}^{P''}(\mathcal{E}))\Big)
\end{equation*}
is a proper closed subscheme of $S$. Let $U_{-}$ be its complement. Let $U_{tf}$ be the dense open subscheme of $S$ consisting of points $s$ such that $\ker\alpha_s$ is a pure sheaf of dimension $d.$ Put $V=U_{-}\cap U_{tf}.$

Suppose that $P_{-}\in C_2$, the other case is similar. By definition of $P_{-}$ the projective morphism $\pi_{fr}\colon \mathrm{FQuot}_{\mathcal{X}/S}^{P_{-}}(\mathcal{E}) \rightarrow S$ is surjective. For any point $s\in S$ the fibre of $\pi_{fr}$ at $s$ parametrizes possible quotients of $E_s$ with Hilbert polynomial $P_{-}.$ If $s\in V$, then any such quotient is a minimal destabilizing framed quotient by construction of $V.$ Recall that the minimal destabilizing framed quotient is unique by Proposition \ref{prop:minquotient}: this implies that $\pi_{fr}\vert_U\colon U:=\pi_{fr}^{-1}(V)\rightarrow V$ is bijective. By the same arguments as in the nonframed case (see \cite[Proposition 3]{art:langton1975}), that quotient is defined over the residue field $k(s)$, hence for $t\in U$, $s=\pi_{fr}(t)$ one has $k(s)\simeq k(t).$ Let $t\in \pi_{fr}^{-1}(s)$ be a point corresponding to a diagram
\begin{equation*}
  \begin{tikzpicture}
    \def\x{1.8}
    \def\y{-1.5}
    \node (A0_0) at (0*\x, 0*\y) {$0$}; 
    \node (A0_1) at (1*\x, 0*\y) {$K$};
    \node (A0_2) at (2*\x, 0*\y) {$E_t$};
    \node (A0_3) at (3*\x, 0*\y) {$Q$};
    \node (A0_4) at (4*\x, 0*\y) {$0$};
    \node (A1_2) at (2*\x, 1*\y) {$F_t$};
    \path (A0_0) edge [->] node [auto] {$\scriptstyle{}$} (A0_1);
    \path (A0_1) edge [->] node [auto] {$\scriptstyle{i}$} (A0_2);
    \path (A0_2) edge [->] node [auto] {$\scriptstyle{q}$} (A0_3);
    \path (A0_3) edge [->] node [auto] {$\scriptstyle{}$} (A0_4);
    \path (A0_2) edge [->] node [auto] {$\scriptstyle{\alpha_t}$} (A1_2);
    \path (A0_3) edge [->] node [auto] {$\scriptstyle{\tilde{\alpha}}$} (A1_2);
  \end{tikzpicture}
\end{equation*}	
By Proposition \ref{prop:tangentefibra}, the Zariski tangent space of $\pi_{fr}^{-1}(s)$ at $t$ is $\mathrm{Hom}(\mathcal{K}, \mathcal{Q})$. Moreover, $K$ is the maximal destabilizing framed subsheaf of $\mathcal{E}_t$, hence $\mathrm{Hom}(\mathcal{K}, \mathcal{Q})=0$ by Proposition \ref{prop:destquo} and therefore $\Omega_{U/V}=0$, hence $\pi_{fr}\vert_U\colon U\rightarrow V$ is unramified. Since $\pi_{fr}$ is projective and $V$ is integral, $\pi_{fr}\vert_U$ is an isomorphism. Now let $T$ be the closure of $U$ in $\mathrm{FQuot}_{\mathcal{X}/S}^{P_{-}}(\mathcal{E})$ with its reduced subscheme structure and $g:=\pi_{fr}\vert_T\colon T\rightarrow S$ is a projective birational morphism. We put $Q$ equal to the pull-back on $X_T$ of the universal quotient on $\mathrm{FQuot}_{\mathcal{X}/S}^{P_{-}}(\mathcal{E})\times_S \mathcal{X}.$

The proof of the universality of the pair $(g,Q)$ is similar to that for the nonframed case (second part of \cite[Theorem 2.3.2]{book:huybrechtslehn2010}), since to prove this part of the theorem we need only the universal property of $\mathrm{FQuot}_{\mathcal{X}/S}^{P_{-}}(\mathcal{E})$ or $\mathrm{Quot}_{\mathcal{X}/S}^{P_{-}}(E).$
\qedhere\endproof
By using the same arguments as in the nonframed case, from the previous theorem one can construct the relative version of the Harder-Narasimhan filtration:
\begin{theorem}
Let $(f\colon \mathcal{X}\rightarrow S, \mathcal{O}_\mathcal{X}(1))$, $F$ and $\delta$ be as in the assumptions of this section. Let $\mathcal{E}=(E, \alpha)$ be a flat family of $d$-dimensional framed sheaves on the fibres of $f$ such that $\ker\alpha_s$ is a pure sheaf of dimension $d$ for a general point $s\in S.$ Assume that if $\beta_d(\mathrm{Im}\:\alpha_s)=0$, then $P(\mathrm{Im}\:\alpha_s)\geq \delta$ for any $s\in S.$ Then there is an integral $k$-scheme $T$ of finite type, a projective birational morphism $g\colon T\rightarrow S$ and a filtration
\begin{equation*}
\mathrm{HN}_{\bullet}(\mathcal{E}): 0=\mathrm{HN}_{0}(\mathcal{E})\subset \mathrm{HN}_{1}(\mathcal{E})\subset \cdots \subset \mathrm{HN}_{l}(\mathcal{E})=E_T
\end{equation*}
such that the following holds:
\begin{itemize}
\item The factors $\mathrm{HN}_{i}(\mathcal{E})/\mathrm{HN}_{i-1}(\mathcal{E})$ are $T$-flat for all $i=1,\ldots, l$, and
\item there is a dense open subscheme $U\subset T$ such that $(\mathrm{HN}_{\bullet}(\mathcal{E}))_t=g_{\mathcal{X}}^*\mathrm{HN}_{\bullet}(\mathcal{E}_{g(t)})$ for all $t\in U.$
\end{itemize}

Moreover, the pair $(g, \mathrm{HN}_{\bullet}(\mathcal{E}))$ is universal in the sense that if $\bar{g}\colon \bar{T}\rightarrow S$ is any dominant morphism of $k$-integral schemes and $\bar{E}_\bullet$ is a filtration of $E_{\bar{T}}$ satisfying these two properties, then there is an $S$-morphism $h\colon \bar{T}\rightarrow T$ such that $h_{\mathcal{X}}^{*}(\mathrm{HN}_{\bullet}(\mathcal{E}))=\bar{E}_\bullet.$
\end{theorem}

\section{$\mu$-(semi)stability}

In this section we give a generalization to framed sheaves of the $\mu$-(semi)stability condition for pure sheaves of dimension $d$ on a projective scheme of dimension $d$ (see \cite[Definition 1.2.12]{book:huybrechtslehn2010}). Also in this case one can construct examples of framed sheaves that are semistable with respect to this new condition but the underlying sheaves are not $\mu$-semistable and vice versa.

Let $(X,\mathcal{O}_X(1))$ be a polarized scheme of dimension $d$ and $\mathcal{E}=(E, \alpha)$ a framed sheaf on it. We define the \textit{rank} (resp.\ the \textit{degree}) of $\mathcal{E}$ as the rank (resp.\ the degree) of $E.$ The \textit{framed degree of} $\mathcal{E}$ is
\begin{equation*}
\deg(\mathcal{E}):=\deg(E)-\epsilon(\alpha)\delta_1.
\end{equation*}
If $E$ has positive rank, its \textit{framed slope} is
\begin{equation*}
\mu(\mathcal{E}):=\frac{\deg(\mathcal{E})}{\mathrm{rk}(\mathcal{E})}.
\end{equation*}
\begin{definition}\label{def:musemi}
A $d$-dimensional framed sheaf $\mathcal{E}=(E,\alpha)$ is $\mu$-(semi)stable with respect to $\delta_1$ if and only if $\ker\alpha$ is a pure sheaf of dimension $d$ and the following conditions are satisfied:
\begin{itemize}
\item[(a)] $\mathrm{rk}(E)\deg(E')\; (\leq)\; \mathrm{rk}(E')\deg(\mathcal{E})$ for all subsheaves $E'\subseteq \ker\alpha$,
\item[(b)] $\mathrm{rk}(E)(\deg(E')-\delta_1)\; (\leq)\; \mathrm{rk}(E')\deg(\mathcal{E})$ for all subsheaves $E'\subset E$ with $\mathrm{rk}(E')< \mathrm{rk}(E).$
\end{itemize}
\end{definition}
One has the usual implications among different stability properties of a framed sheaf of positive rank:
\begin{equation*}
\mu\mbox{-stable}\Rightarrow \mbox{stable} \Rightarrow \mbox{semistable}\Rightarrow \mu\mbox{-semistable}.
\end{equation*}
\begin{definition}
Let $\mathcal{E}=(E,\alpha)$ be a framed sheaf of dimension $d-1.$ If $\alpha$ is injective, we say that $\mathcal{E}$ is $\mu$-semistable\footnote{For sheaves of dimension $d-1$, the definition of $\mu$-semistability of the corresponding framed sheaves does not depend on $\delta_1.$}. Moreover, if the degree of $E$ is $\delta_1$, we say that $\mathcal{E}$ is $\mu$-stable with respect to $\delta_1.$
\end{definition}
\begin{remark}\label{rem:results}
Most results of Sections \ref{sec:semistability} - \ref{sec:relative} still hold if one replaces (semi)stability by $\mu$-(semi)stability. For example, this is the case for Theorem \ref{thm:rel}, Propositions \ref{prop:quoziente1}, \ref{prop:dest}, \ref{prop:minquotient}, \ref{prop:extsocle}, Corollaries \ref{cor:semist}, \ref{cor:geostable}. 
\end{remark}
Below we state separately the theorem on the existence of the $\mu$-Harder-Narasimhan and $\mu$-Jordan-H\"older filtrations. We start by specifying the definitions of these notions.
\begin{definition}
$\mbox{ }$
\begin{itemize}
 \item Let $\mathcal{E}=(E, \alpha)$ be a framed sheaf where $\ker\alpha$ is torsion free and $\mathrm{Im}\:\alpha$ has rank zero. A $\mu$\textit{-Harder-Narasimhan filtration} is a filtration of framed saturated subsheaves
\begin{equation*}
\mathrm{HN}_{\bullet}(\mathcal{E}): 0=\mathrm{HN}_{0}(\mathcal{E})\subset \mathrm{HN}_{1}(\mathcal{E})\subset \cdots \subset \mathrm{HN}_{l}(\mathcal{E})=E
\end{equation*}
which satisfies the following conditions
\begin{itemize}
\item the quotient sheaf $gr_i^{\mathrm{HN}}(\mathcal{E}):=\mathrm{HN}_{i}(\mathcal{E})/\mathrm{HN}_{i-1}(\mathcal{E})$ with the induced framing $\alpha_i$ is a $\mu$-semistable framed sheaf for $i=1, 2, \ldots, l.$
\item the quotient $\left(E/\mathrm{HN}_{j}(\mathcal{E}), \alpha''\right)$ is a framed sheaf where $\ker\alpha''$ is torsion free for $j=1,2, \ldots, l-2$, it has no rank zero $\mu$-destabilizing framed subsheaves, and
\begin{equation*}
\mathrm{rk}(gr_{i+1}^{\mathrm{HN}}(\mathcal{E}))\deg(gr_i^{\mathrm{HN}}(\mathcal{E}), \alpha_i)>\mathrm{rk}(gr_{i}^{\mathrm{HN}}(\mathcal{E}))\deg(gr_{i+1}^{\mathrm{HN}}(\mathcal{E}), \alpha_{i+1})
\end{equation*}
for $i=1, \ldots, l-1.$
\end{itemize}
\item Let $\mathcal{E}=(E,\alpha)$ be a $\mu$-semistable framed sheaf of positive rank. A $\mu$\textit{-Jordan-H\"older filtration} is a filtration of framed saturated subsheaves
\begin{equation*}
E_{\bullet}: 0=E_0\subset E_1\subset\cdots \subset E_l=E
\end{equation*}
such that all the factors $E_i/E_{i-1}$ together with the induced framings $\alpha_i$ are $\mu$-stable with framed degree $\deg(E_i/E_{i-1},\alpha_i)=\mathrm{rk}(E_i/E_{i-1})\mu(\mathcal{E}).$ 
\end{itemize}
\end{definition}
\begin{theorem}
$\mbox{ }$
\begin{itemize}
 \item Let $\mathcal{E}=(E, \alpha)$ be a framed sheaf where $\ker\alpha$ is torsion free and $\mathrm{Im}\:\alpha$ has rank zero. Then there exists a unique $\mu$-Harder-Narasimhan filtration.
\item Let $\mathcal{E}=(E,\alpha)$ be a $\mu$-semistable framed sheaf of positive rank. Then there exist $\mu$-Jordan-H\"older filtrations. Moreover, the graded framed sheaf
\begin{equation}\label{eq:graded}
gr(\mathcal{E})=(gr(E),gr(\alpha))=\bigoplus_i (E_i/E_{i-1},\alpha_i)
\end{equation}
does not depend on the choice of the $\mu$-Jordan-H\"older filtration.
\end{itemize}
\end{theorem}
For $i\geq 0$, let us denote by $\mathrm{Coh}_{i}(X)$ the full subcategory of $\mathrm{Coh}(X)$ whose objects are sheaves of dimension less or equal to $i.$ 

Let $\mathrm{Coh}_{d,d-1}(X)$ be the quotient category $\mathrm{Coh}_{d}(X)/\mathrm{Coh}_{d-1}(X).$ In \cite[Section 1.6]{book:huybrechtslehn2010}, Huybrechts and Lehn define the notion of $\mu$-Jordan-H\"older filtration for $\mu$-semistable sheaves $E$ in the category $\mathrm{Coh}_{d,d-1}(X).$ For a $\mu$-semistable torsion free sheaf $E$, the graded object associated to a $\mu$-Jordan-H\"older filtration is uniquely determined only in codimension one. In our case, we construct $\mu$-Jordan-H\"older filtrations by using filtrations in which every term is a framed saturated subsheaf of the next term. In this way, the graded object \eqref{eq:graded} is uniquely determined.  Thus, when the framing of a $\mu$-semistable framed sheaf is zero, our definition of $\mu$-Jordan-H\"older filtration does not coincide with the nonframed one given by Huybrechts and Lehn.

\section{Restriction theorems}\label{sec:restriction}

In this section we generalize the Mehta-Ramanathan restriction theorems to framed sheaves. We limit our attention to the case in which the framing sheaf $F$ is a sheaf supported on a divisor $\FD.$ In the framed case the results depend also on the parameter $\delta_1.$ We shall follow rather closely the arguments and the techniques used by Huybrechts and Lehn to prove the Mehta-Ramanathan restriction theorems in the nonframed case (see \cite[Section 7.2]{book:huybrechtslehn2010}). In the framed case, the proofs are somewhat more elaborate than in the nonframed case.

\subsubsection*{$\mu$-semistable case}

We provide a generalization of Mehta-Ramanathan's restriction theorem for $\mu$-semista\-ble torsion free sheaves \cite[Theorem 6.1]{art:mehtaramanathan1981}.
\begin{theorem}\label{thm:mr1}
Let $(X,\mathcal{O}_X(1))$ be a nonsingular polarized variety of dimension $d\geq 2.$ Let $F$ be a sheaf on $X$ supported on a divisor $\FD.$ Let $\mathcal{E}=(E,\alpha\colon E\rightarrow F)$ be a framed sheaf on $X$ of positive rank with nonzero framing. If $\mathcal{E}$ is $\mu$-semistable with respect to $\delta_1$, there exists a positive integer $a_0$ such that for all $a\geq a_0$ there is a dense open subset $U_a\subset \vert \mathcal{O}_X(a)\vert$ such that for all $D\in U_a$ the divisor $D$ is smooth, meets transversely the divisor $\FD$ and $\mathcal{E}\vert_D$ is $\mu$-semistable with respect to $a\delta_1.$
\end{theorem}
\proof
Suppose the theorem is false. Then for any positive integer $a$, there exists a general divisor $D\in \vert \mathcal{O}_X(a)\vert$ such that $\mathcal{E}\vert_D$ is not $\mu$-semistable with respect to $a\delta_1.$ This is equivalent to say that there exists a framed quotient $\mathcal{E}\vert_D\rightarrow \mathcal{G}$ such that $\mu(\mathcal{G})<\mu(\mathcal{E}\vert_D).$ Our idea is to extend this quotient to a $\mu$-destabilizing framed quotient of $\mathcal{E}.$ The main tool we shall use to do this is Theorem \ref{thm:rel}, which states the existence of the relative minimal $\mu$-destabilizing framed quotient with respect to $a\delta_1.$ We shall divide the proof in several steps.

\emph{Step 1:} first, we define the setting in which we use Theorem \ref{thm:rel}. 

For a positive integer $a$, let $\Pi_a:=\vert \mathcal{O}_X(a)\vert$ be the complete linear system of hypersurfaces of degree $a$ and let $Z_a:=\{(D,x)\in \vert \mathcal{O}_X(a)\vert\times X\mid x\in D\}$ be the \textit{incidence variety} with its natural projections
\begin{equation*}
  \begin{tikzpicture}
    \def\x{1.5}
    \def\y{-1.2}
    \node (A0_0) at (0*\x, 0*\y) {$Z_a$};
    \node (A0_1) at (1*\x, 0*\y) {$X$};
    \node (A1_0) at (0*\x, 1*\y) {$\Pi_a$};
    \path (A0_0) edge [->] node [auto] {$\scriptstyle{q}$} (A0_1);
    \path (A0_0) edge [->] node [auto] {$\scriptstyle{p}$} (A1_0);
  \end{tikzpicture}
\end{equation*}
It is possible to give a schematic structure on $Z_a$ so that $p$ is a projective morphism with equidimensional fibres (see \cite[Section 3.1]{book:huybrechtslehn2010}). Thus we can induce on $p\colon Z_a\rightarrow \Pi_a$ a structure of relative projective scheme of fibre dimension $d-1.$ Moreover one can prove the following property (see \cite[Section 2]{art:mehtaramanathan1981}, \cite[Exercise II.6.1]{book:hart}):
\begin{equation}\label{eq:decomposition}
\mathrm{Pic}(Z_a)=q^*(\mathrm{Pic}(X))\oplus p^*(\mathrm{Pic}(\Pi_a)).  
\end{equation}
\begin{lemma}
The pair $(q^{*}E, q^*\alpha)$ is a flat family of $(d-1)$-dimensional framed sheaves on the fibres of $p.$ 
\end{lemma}
\proof
For all $D\in \Pi_a$, the Hilbert polynomials of the restrictions $E\vert_D$, $F\vert_D$ and $\mathrm{Im}\:\alpha\vert_D$ are independent of $D$, indeed, e.g., the Hilbert polynomial of $E\vert_D$ is $P(E\vert_D,n)=P(E,n)-P(E,n-a).$ Since $\Pi_a$ is a reduced scheme, by \cite[Proposition 2.1.2]{book:huybrechtslehn2010} $q^{*}F$, $q^*E$ and $q^*\mathrm{Im}\:\alpha$ are flat families of sheaves on the fibres of $p.$ Moreover, for all $D\in \Pi_a$ we have $\alpha\vert_D\neq 0$, $\dim(E\vert_D)=d-1$ and  $\dim(F\vert_D)<d-1.$
\black

\emph{Step 2:} we apply Theorem \ref{thm:rel} to the flat family $(q^{*}E, q^*\alpha).$

For general $D\in \Pi_a$ the restriction $\ker\alpha\vert_D$ is again torsion free (see \cite[Corollary 1.1.14]{book:huybrechtslehn2010}), hence the open subset of $\Pi_a$ consisting of divisors $D$ such that $\ker\alpha\vert_D$ is pure of dimension $d-1$ is nonempty. Since $\mathcal{E}$ is $\mu$-semistable with respect to $\delta_1$, $\deg(\mathrm{Im}\:\alpha)\geq \delta_1$, hence $\deg(\mathrm{Im}\:\alpha\vert_D)=a\deg(\mathrm{Im}\:\alpha)\geq a\delta_1.$

Thus the flat family $(q^{*}E, q^*\alpha)$ satisfies the hypothesis of Theorem \ref{thm:rel}. Hence there are a dense open subset $V_a\subset \Pi_a$ and a $V_a$-flat quotient on $Z_{V_a}:=p^{-1}(V_a)=V_a\times_{\Pi_a} Z_a$
\begin{equation}\label{eq:quotient}
\begin{aligned}
  \begin{tikzpicture}
    \def\x{2}
    \def\y{-1.8}
    \node (A0_0) at (0*\x, 0*\y) {$(q^{*}E)\vert_{Z_{V_a}}$};
    \node (A0_1) at (1*\x, 0*\y) {$Q_a$};
    \node (A1_0) at (0*\x, 1*\y) {$(q^*F)\vert_{Z_{V_a}}$};
    \path (A0_0) edge [->] node [auto] {$q_a$} (A0_1);
    \path (A0_0) edge [->] node [left] {$(q^*\alpha)\vert_{Z_{V_a}}$} (A1_0);
    \path (A0_1) edge [->] node [auto] {$\tilde{\alpha}_a$} (A1_0);
  \end{tikzpicture} 
\end{aligned}
\end{equation}
with a morphism $\tilde{\alpha}_a\colon Q_a\rightarrow (q^*F)\vert_{Z_{V_a}}$, such that for all $D\in V_a$ the framed sheaf $(E\vert_D,\alpha\vert_D)$ has dimension $d-1$, and $\ker\alpha\vert_D$ is a pure sheaf of dimension $d-1$; moreover, $Q_a\vert_D$ is a $(d-1)$-dimensional sheaf, $\tilde{\alpha}_a\vert_D$ is the framing induced by $\alpha\vert_D$ and $(Q_a\vert_D,\tilde{\alpha}_a\vert_D)$ is the minimal $\mu$-destabilizing framed quotient of $(E\vert_D, \alpha\vert_D).$ 

\emph{Step 3:} we need to introduce the quantities $\mathrm{rk}(a)$ and $\mu_{fr}(a).$

Fix an extension of $\det(Q_a)$ to some line bundle on all of $Z_a.$ This can be uniquely decomposed in the form $q^*L_a\otimes p^*M$ with $L_a\in \mathrm{Pic}(X)$ and $M\in \mathrm{Pic}(\Pi_a).$ Note that $\deg(Q_a\vert_D)=a\deg(L_a)$ for $D\in V_a.$ 

Let $U_a\subset V_a$ be the dense open set of points $D\in V_a$ such that $D$ is smooth and $D$ meets transversely the divisor $\FD.$ Let $\deg(a)$, $\mathrm{rk}(a)$ and $\mu_{fr}(a)$ denote the degree, the rank and the framed slope of the minimal $\mu$-destabilizing framed quotient of $(E\vert_D, \alpha\vert_D)$ for a general point $D\in \Pi_a.$ By construction of the relative minimal $\mu$-destabilizing framed quotient, the quantity $\epsilon(\tilde{\alpha}_a\vert_D)$ is independent of $D\in V_a$, so we denote it by $\epsilon(a).$ Then we have $1\leq \mathrm{rk}(a)\leq \mathrm{rk}(E)$ and
\begin{equation*}
\frac{\mu_{fr}(a)}{a}=\frac{\deg(Q_a\vert_D)-\epsilon(a)a\delta_1}{\mathrm{rk}(a)a}=\frac{\deg(L_a)-\epsilon(a)\delta_1}{\mathrm{rk}(a)}\in\frac{\mathbb{Z}}{\delta_1''(\mathrm{rk}(E)!)}\subset \mathbb{Q},
\end{equation*}
where $\delta_1=\delta_1'/\delta_1''.$ Let $l>1$ be an integer, $a_1,\ldots, a_l$ positive integers and $a=\sum_i a_i.$ 

\emph{Step 4:} we prove that $\mathrm{rk}(a)$ and $\mu_{fr}(a)/a$ are eventually constant.

First, we need to compare $\mathrm{rk}(a)$ (resp.\ $\mu_{fr}(a)/a$) with $\mathrm{rk}(a_i)$ (resp.\ $\mu_{fr}(a_i)/a_i$) for all $i=1, \ldots, l.$ To do this, we use the following result, which allows us to compare the rank and the framed degree of $Q_{a_i}$ in a generic fibre with the same invariants of a ``special quotient'' of $(q^{*}E)\vert_{Z_{V_a}}.$
\begin{lemma}{\normalfont (\cite[Lemma 7.2.3]{book:huybrechtslehn2010}).}\label{lem:divi}
Let $l>1$ be an integer, $a_1, \ldots, a_l$ positive integers, $a=\sum_i a_i$, and let $D_i\in U_{a_i}$ be divisors such that $D=\sum_i D_i$ is a divisor with normal crossings. Then there is a smooth locally closed curve $C\subset \Pi_a$ containing the point $D\in \Pi_a$ such that $C\setminus \{D\}\subset U_a$ and $Z_C:=p^{-1}(C)=C\times_{\Pi_a} Z_a$ is smooth in codimension 2.
\end{lemma}
Note that if $D_1\in U_{a_1}$ is given, one can always find $D_i\in U_{a_i}$ for $i\geq 2$ such that $D=\sum_i D_i$ is a divisor with normal crossings.
\begin{lemma}\label{lem:sum}
Let $a_1,\ldots, a_l$ be positive integers, with $l>1$, and $a=\sum_i a_i.$ Then $\mu_{fr}(a)\geq \sum_i \mu_{fr}(a_i)$ and in case of equality $\mathrm{rk}(a)\geq \max\{\mathrm{rk}(a_i)\}.$
\end{lemma}
\proof
Let $D_i$ be divisors satisfying the requirements of Lemma \ref{lem:divi} and let $C$ be the curve with the properties of Lemma \ref{lem:divi}. 

Now we have to consider two cases for the quotient $(q^{*}E)\vert_{Z_{V_a}}\stackrel{q_a}{\longrightarrow} Q_a$ introduced in diagram \eqref{eq:quotient}:
\begin{itemize}
\item[(1)] there exists a framing $\tilde{\alpha}_a$ on $Q_a$ such that $(q^*\alpha)\vert_{Z_{V_a}}=\tilde{\alpha}_a\circ q_a$,
\item[(2)] $\ker q_a\vert_{D'}\not\subset\ker\alpha\vert_{D'}$ for all $D'\in V_a.$
\end{itemize}

In the first case, $\tilde{\alpha}_a\vert_{D'}\neq 0$ for all $D'\in V_a.$ The restriction of diagram \eqref{eq:quotient} to $Z_{V_a\cap C}$ is
\begin{equation*}
  \begin{tikzpicture}
    \def\x{3.4}
    \def\y{-2.1}
    \node (A0_0) at (0*\x, 0*\y) {$0$};
    \node (A0_1) at (1*\x, 0*\y) {$K$};
    \node (A0_2) at (2*\x, 0*\y) {$(q^{*}E)\vert_{Z_{V_a\cap C}}$};
    \node (A0_3) at (3*\x, 0*\y) {$Q_a\vert_{Z_{V_a\cap C}}$};
    \node (A0_4) at (4*\x, 0*\y) {$0$};
    \node (A2_1) at (2*\x, 1*\y) {$(q^*F)\vert_{Z_{V_a\cap C}}$};
     \path (A0_0) edge [->] node [auto] {$$} (A0_1);
      \path (A0_1) edge [->] node [auto] {$$} (A0_2);
    \path (A0_2) edge [->] node [auto] {$q_a\vert_{Z_{V_a\cap C}}$} (A0_3);
     \path (A0_3) edge [->] node [auto] {$$} (A0_4);
    \path (A0_2) edge [->] node [left=2pt] {$(q^*\alpha)\vert_{Z_{V_a\cap C}}$} (A2_1);
    \path (A0_3) edge [->] node [auto] {$\tilde{\alpha}_a\vert_{Z_{V_a\cap C}}$} (A2_1);
  \end{tikzpicture}
\end{equation*}

Since the morphism $Z_{V_a\cap C}\rightarrow Z_C$ is flat (because it is an open embedding), $\ker (q^*\alpha\vert_{Z_{V_a\cap C}})=(\ker q^{*}\alpha\vert_{Z_C})\vert_{Z_{V_a\cap C}}$ and we can extend the inclusion $K\subset \ker q^*\alpha\vert_{Z_{V_a\cap C}}$ to an inclusion $K_C\subset \ker q^{*}\alpha\vert_{Z_C}$ on $Z_C.$ Therefore we can extend $Q_a\vert_{Z_{V_a\cap C}}$ to a $C$-flat quotient $Q_C$ of $q^{*}E\vert_{Z_C}$ such that there exists a framing $\tilde{\alpha}_C$ on $Q_C$ for which the following diagram commutes
\begin{equation*}
  \begin{tikzpicture}
    \def\x{2.3}
    \def\y{-2.1}
    \node (A0_0) at (0*\x, 0*\y) {$(q^{*}E)\vert_{Z_{C}}$};
    \node (A0_1) at (1*\x, 0*\y) {$Q_C$};
    \node (A1_0) at (0*\x, 1*\y) {$(q^*F)\vert_{Z_{C}}$};
    \path (A0_0) edge [->] node [auto] {$q_C$} (A0_1);
    \path (A0_0) edge [->] node [left=2pt] {$(q^*\alpha)\vert_{Z_{C}}$} (A1_0);
    \path (A0_1) edge [->] node [auto] {$\tilde{\alpha}_C$} (A1_0);
  \end{tikzpicture}
\end{equation*}
Furthermore $\tilde{\alpha}_C\vert_c\neq 0$ for all $c\in C.$ By the flatness of $Q_C$ we obtain $P(Q_C\vert_c, n)=P(Q_C\vert_D,n)$ for all $c\in C\setminus \{D\}.$ Hence $\mathrm{rk}(Q_C\vert_D)=\mathrm{rk}(a)$ and $\deg(Q_C\vert_D)=\deg(a)$, therefore $\mu(Q_C\vert_D,\tilde{\alpha}_C\vert_D)=\mu_{fr}(a).$ 

Let $T'(Q_C\vert_D)$ be the sheaf on $D$ that to every open subset $U$ associates the set of sections $f$ of $Q_C\vert_D$ in $U$ such that there exists $n>0$ for which $I_D^n\cdot f=0$, where $I_D$ is the ideal sheaf associated to $D.$ Roughly speaking, $T'(Q_C\vert_D)$ is the part of the torsion subsheaf $T(Q_C\vert_D)$ of $Q_C\vert_D$ that is not supported in the intersection $D\cap \FD.$ By the transversality of $D_i$ with respect to $\FD$, we have $T'(Q_C\vert_D)\subset \ker\tilde{\alpha}_C\vert_D.$ Thus the sheaf $\bar{Q}=Q_C\vert_D/T'(Q_C\vert_D)$ is a positive rank quotient of $E\vert_D$ with nonzero induced framing $\bar{\alpha}.$ Moreover, $\mathrm{rk}(\bar{Q}\vert_{D_i})=\mathrm{rk}(\bar{Q})=\mathrm{rk}(Q_C\vert_D)=\mathrm{rk}(a).$ Since $T'(Q_C\vert_D)$ is a rank zero sheaf, the $(d-2)$-th coefficient $\beta_{d-2}(T'(Q_C\vert_D))$ of its Hilbert polynomial is nonnegative, hence $\deg(Q_C\vert_D)=\deg(\bar{Q})+\beta_{d-2}(T'(Q_C\vert_D))\geq \deg(\bar{Q}).$ The framing of $\bar{Q}$ is nonzero, hence the previous inequality between degrees yields an inequali\-ty between framed degrees. Therefore
\begin{equation*}
\mu_{fr}(a)=\mu(Q_C\vert_D,\tilde{\alpha}_C\vert_D)\geq \mu(\bar{Q},\bar{\alpha}).
\end{equation*}
The sequence
\begin{equation*}
0\longrightarrow \bar{Q}\longrightarrow \bigoplus_i \bar{Q}\vert_{D_i}\longrightarrow \bigoplus_i\bigoplus_{i<j}\bar{Q}\vert_{D_i\cap D_j}\longrightarrow 0
\end{equation*}
is exact modulo sheaves of dimension $d-3$ (the kernel of the morphism $\bar{Q}\longrightarrow \bigoplus_i \bar{Q}\vert_{D_i}$ is zero because the divisors $D_i$ are transversal with respect to the singular set of $\bar{Q}$). By the same computations as in the proof of \cite[Lemma 7.2.5]{book:huybrechtslehn2010} we have
\begin{equation*}
\mu(\bar{Q})=\sum_i\Big(\mu(\bar{Q}\vert_{D_i})-\frac{1}{2}\sum_{j\neq i}\Big(\frac{\mathrm{rk}\left(\bar{Q}\vert_{D_i\cap D_j}\right)}{\mathrm{rk}(a)}-1\Big)a_i a_j\Big).
\end{equation*}
For every $i$ and $j\neq i$ we also define the sheaf $T_{ij}(\bar{Q}\vert_{D_i})$ as the sheaf on $D_i$ that to every open subset $U$ associates the set of sections $f$ of $\bar{Q}\vert_{D_i}$ in $U$ such that there exists $n>0$ for which $I_{D_j}^n\cdot f=0.$ Note that $T_{ij}(\bar{Q}\vert_{D_i})\subset \ker\bar{\alpha}\vert_{D_i}.$ We define $Q_i=\bar{Q}\vert_{D_i}/\bigoplus_{j\neq i}T_{ij}(\bar{Q}\vert_{D_i}).$ By construction $\mathrm{rk}(Q_i)=\mathrm{rk}(\bar{Q})$ and there exists a nonzero induced framing $\alpha_i$ on $Q_i.$ As before, by the same computations as in the proof of \cite[Lemma 7.2.5]{book:huybrechtslehn2010}, we obtain
\begin{equation*}
\mu(Q_i)=\mu(\bar{Q}\vert_{D_i})-\sum_{j\neq i}\Big(\frac{\mathrm{rk}\left(\bar{Q}\vert_{D_i\cap D_j}\right)}{\mathrm{rk}(a)}-1\Big)a_i a_j.
\end{equation*}
Therefore $\mu(\bar{Q})\geq \sum_i \mu(Q_i)$, and 
\begin{equation*}
\mu_{fr}(a)\geq \mu(\bar{Q},\bar{\alpha})\geq \sum_i \mu(Q_i,\alpha_i).
\end{equation*}
By definition of minimal $\mu$-destabilizing framed quotient, we have $\mu(Q_i,\alpha_i)\geq \mu_{fr}(a_i)$, hence $\mu_{fr}(a)\geq \sum_i\mu_{fr}(a_i).$

Consider the second case. On the restriction to $Z_{V_a\cap C}$ there is the quotient:
\begin{equation*}
  \begin{tikzpicture}
    \def\x{3.5}
    \def\y{-1.8}
    \node (A0_0) at (0*\x, 0*\y) {$(q^{*}E)\vert_{Z_{V_a\cap C}}$};
    \node (A0_1) at (1*\x, 0*\y) {$Q_a\vert_{Z_{V_a\cap C}}$};
    \node (A1_0) at (0*\x, 1*\y) {$(q^*F)\vert_{Z_{V_a\cap C}}$};
    \path (A0_0) edge [->] node [auto] {$q$} (A0_1);
    \path (A0_0) edge [->] node [auto] {$(q^*\alpha)\vert_{Z_{V_a\cap C}}$} (A1_0);
  \end{tikzpicture}
\end{equation*}
By definition of $Q_a$, we get $\ker q\vert_{D'}\not\subset \ker \alpha\vert_{D'}$ for all points $D'\in V_a\cap C$, hence $\ker q\not\subset \ker (q^*\alpha)\vert_{Z_{V_a\cap C}}.$ As before, we can extend $Q_a\vert_{Z_{V_a\cap C}}$ to a $C$-flat quotient
\begin{equation*}
  \begin{tikzpicture}
    \def\x{2}
    \def\y{-1.8}
    \node (A0_0) at (0*\x, 0*\y) {$(q^{*}E)\vert_{Z_{C}}$};
    \node (A0_1) at (1*\x, 0*\y) {$Q_C$};
    \node (A1_0) at (0*\x, 1*\y) {$(q^*F)\vert_{Z_{C}}$};
    \path (A0_0) edge [->] node [auto] {$q_C$} (A0_1);
    \path (A0_0) edge [->] node [auto] {$(q^*\alpha)\vert_{Z_{C}}$} (A1_0);
  \end{tikzpicture}
\end{equation*}
Since $\ker q_C$ and $\ker (q^*\alpha)\vert_{Z_{C}}$ are $C$-flat, also $\ker q_C\cap \ker (q^*\alpha)\vert_{Z_{C}}$ is $C$-flat. Moreover for all points $D'\in V_a\cap C$ we have $(\ker q_C\cap \ker(q^*\alpha)\vert_{Z_{C}})\vert_{D'}=\ker q\vert_{D'}\cap \ker \alpha\vert_{D'}$, hence by flatness we get $\ker q_C\vert_{D'}\not\subset \ker \alpha\vert_{D'}$ for all points $D'\in C.$ As before, by flatness of $Q_C$ we have that $\mathrm{rk}(Q_C\vert_D)=\mathrm{rk}(a)$ and $\deg(Q_C\vert_D)=\deg(a)$; moreover the induced framing on $Q_C\vert_D$ is zero, hence $\mu(Q_C\vert_D)=\mu_{fr}(a).$ Let $\bar{Q}=Q_C\vert_D/T(Q_C\vert_D)$ and $Q_i=\bar{Q}\vert_{D_i}/T(\bar{Q}\vert_{D_i}).$ Using the same computations as in the proof of \cite[Lemma 7.2.5]{book:huybrechtslehn2010}, we obtain $\mu(\bar{Q})\geq \sum_i\mu(Q_i).$ As before, we get $\mu_{fr}(a)=\mu(Q_C\vert_D)\geq \mu(\bar{Q})\geq \sum_i \mu(Q_i)\geq \sum_i \mu_{fr}(a_i).$

Now let us consider the case $\mu_{fr}(a)=\sum_i \mu_{fr}(a_i).$ In both cases, if we denote by $\alpha_i$ the induced framing on $Q_i$, from this equality, follows that $\mu(Q_i, \alpha_i)=\mu_{fr}(a_i)$ and $\mathrm{rk}(\bar{Q}\vert_{D_i\cap D_j})=\mathrm{rk}(a).$ Since $\mu_{fr}(a_i)$ is the framed-slope of the minimal $\mu$-destabilizing framed quotient, we have that $\mathrm{rk}(a)=\mathrm{rk}(Q_i)\geq \mathrm{rk}(a_i)$ for all $i.$
\black

\begin{corollary}\label{cor:1}
$\mathrm{rk}(a)$ and $\mu_{fr}(a)/a$ are constant for $a\gg 0.$
\end{corollary}
\proof
The function $a\mapsto \mu_{fr}(a)/a$ takes values in a lattice of $\mathbb{Q}$ and is bounded from above by $\mu(\mathcal{E}).$ Thus one can use the same arguments as in the nonframed case (Corollary 7.2.6 in \cite{book:huybrechtslehn2010}) and obtain the assertion. In a similar way, one can prove that $\mathrm{rk}(a)$ is eventually constant.
\black

\emph{Step 5:} we prove that the line bundle $L_a$ is eventually constant.

If $\mu_{fr}(a)/a=\mu_{fr}(a_i)/a_i$ and $\mathrm{rk}(a)=\mathrm{rk}(a_i)$ for all $i$, then $Q_i$ is the minimal $\mu$-destabilizing framed quotient of $E\vert_{D_i}$, hence $Q_C\vert_{D_i}$ differs from the minimal $\mu$-destabilizing framed quotient of $E\vert_{D_i}$ only in dimension $d-3$, in particular their determinant line bundles as sheaves on $D_i$ are equal. From this argument and the same techniques as in the proof of the corresponding result for the nonframed case (\cite[Lemma 7.2.7]{book:huybrechtslehn2010}), it follows:
\begin{lemma}\label{lem:2}
There is a line bundle $L\in \mathrm{Pic}(X)$ such that $L_a \simeq L$ for all $a\gg 0.$
\end{lemma}
By the previous lemma and Corollary \ref{cor:1}, $\epsilon(a)$ is constant for $a\gg 0.$   

\emph{Step 6:} we construct a $\mu$-destabilizing framed quotient of $\mathcal{E}.$

Summarizing what we have obtained until now, we proved that the quantity
\begin{equation*}
\frac{\deg(Q_a\vert_D)-\epsilon(\tilde{\alpha}_a\vert_D)a\delta_1}{\mathrm{rk}(a)}=a\frac{\deg(L)-\epsilon(a)\delta_1}{\mathrm{rk}(a)}
\end{equation*}
is independent of $D\in V_a$ and $a\gg 0.$

Now we have all the ingredients to construct a positive rank $\mu$-destabilizing quotient of $\mathcal{E}$, as we shall explain in the following. We have to consider separately two cases: $\epsilon(a)=1$ and $\epsilon(a)=0$ for $a\gg 0.$ In the first case we have
\begin{equation*}
\frac{\deg(L)-\delta_1}{r}<\mu(\mathcal{E})
\end{equation*}
and $1\leq r\leq\mathrm{rk}(E)$, where $r=\mathrm{rk}(a)$ for $a\gg 0.$ We want to construct a rank $r$ quotient $Q$ of $E$, with nonzero induced framing $\beta$ and determinant line bundle equal to the line bundle $L$ defined in Lemma \ref{lem:2}. Then we would get $\mu(\mathcal{Q})<\mu(\mathcal{E})$ and therefore we would obtain a contradiction with the hypothesis of $\mu$-semistability of $\mathcal{E}$ with respect to $\delta_1.$ 

Let $a$ be sufficiently large, $D\in U_a$ and the minimal $\mu$-destabilizing framed quotient
\begin{equation*}
  \begin{tikzpicture}
    \def\x{2}
    \def\y{-1.8}
    \node (A0_0) at (0*\x, 0*\y) {$E\vert_{D}$};
    \node (A0_1) at (1*\x, 0*\y) {$Q_D$};
    \node (A1_0) at (0*\x, 1*\y) {$F\vert_{D}$};
    \path (A0_0) edge [->] node [auto] {$q_D$} (A0_1);
    \path (A0_0) edge [->] node [left=2pt] {$\alpha\vert_{D}$} (A1_0);
    \path (A0_1) edge [->] node [auto] {$\beta_D$} (A1_0);
  \end{tikzpicture}
\end{equation*}
Put $K_D=\ker\beta_D.$ By Proposition \ref{prop:quoziente1} (for $\mu$-semistability), $Q_D$ fits into an exact sequence
\begin{equation}\label{eq:overD}
0\longrightarrow K_D\longrightarrow Q_D\longrightarrow \mathrm{Im}\:\alpha\vert_D\longrightarrow 0
\end{equation}
with $K_D$ torsion free quotient of $\ker\alpha\vert_D.$ So there exists an open subscheme $\tilde{D}\subset D$ such that $K_D\vert_{\tilde{D}}$ is locally free of rank $r$ and $D\setminus \tilde{D}$ is a closed subset of codimension at least two in $D.$ Consider the restriction of the sequence \eqref{eq:overD} on $\tilde{D}$
\begin{equation*}
0\longrightarrow K_D\vert_{\tilde{D}}\longrightarrow Q_D\vert_{\tilde{D}}\longrightarrow \mathrm{Im}\:\alpha\vert_{\tilde{D}}\longrightarrow 0.
\end{equation*}
By taking determinants of the exact sequence, we have
\begin{equation*}
\det(K_D)\vert_{\tilde{D}}\otimes \det(\mathrm{Im}\:\alpha\vert_{\tilde{D}})=\det(Q_D\vert_{\tilde{D}}).
\end{equation*}
Let us denote by $\bar{L}$ the line bundle $L\otimes \det(\mathrm{Im}\:\alpha)^\vee.$ So we get $\det(K_D)\vert_{\tilde{D}}=\bar{L}\vert_{\tilde{D}}.$ Therefore $\det(K_D)=\bar{L}\vert_{D}.$ 

We have a morphism $\sigma_D\colon \Lambda^r\ker\alpha\vert_D\rightarrow \bar{L}\vert_{D}$ which is surjective on $\tilde{D}$ and morphisms
\begin{equation*}
\tilde{D}\longrightarrow \mathrm{Grass}(\ker\alpha,r)\longrightarrow \mathbb{P}(\Lambda^r\ker\alpha).
\end{equation*}
By Serre's vanishing theorem and Serre duality, one has for $i=0,1$
\begin{equation*}
\mathrm{Ext}^i(\Lambda^r\ker\alpha,\bar{L}(-a))=\mathrm{H}^{d-i}(X, \Lambda^r\ker\alpha\otimes\bar{L}^\vee\otimes\omega_X(a))^\vee=0
\end{equation*}
for all $a\gg 0$ (since $d\geq 2$), hence
\begin{equation*}
\mathrm{Hom}(\Lambda^r\ker\alpha,\bar{L})=\mathrm{Hom}(\Lambda^r\ker\alpha\vert_D,\bar{L}\vert_D).
\end{equation*}
So for $a$ sufficiently large, the morphism $\sigma_D$ extends to a morphism $\sigma\colon \Lambda^r\ker\alpha\rightarrow \bar{L}.$ The support of the cokernel of $\sigma$ meets $D$ in a closed subscheme of codimension two in $D$, hence there is an open subscheme $\tilde{X}\subset X$ such that $\sigma\vert_{\tilde{X}}$ is surjective, $X\setminus \tilde{X}$ is a closed subscheme of codimension two and $\tilde{D}=\tilde{X}\cap D.$ So we have a morphism $\tilde{X}\rightarrow \mathbb{P}(\Lambda^r\ker\alpha)$ and we want that it factorizes through $\mathrm{Grass}(\ker\alpha,r).$ Using the same arguments of the final part of proof of \cite[Theorem 7.2.1]{book:huybrechtslehn2010}, for sufficiently large $a$ that morphism factorizes, hence we get a rank $r$ locally free quotient $\ker\alpha\vert_{\tilde{X}}\longrightarrow K_{\tilde{X}}$ such that $\det(K_{\tilde{X}})=\bar{L}\vert_{\tilde{X}}.$ So we can extend $K_{\tilde{X}}$ to a rank $r$ quotient $K$ of $\ker\alpha$ such that $\det(K)=\bar{L}.$ 

Let $G=\ker(\ker\alpha\rightarrow K)$ and $Q=E/G.$ Since $Q$ fits into the exact sequence
\begin{equation*}
 0\longrightarrow K\longrightarrow Q\longrightarrow \mathrm{Im}\:\alpha\longrightarrow 0,
\end{equation*}
the determinant of $Q$ is $\det(K)\otimes\det(\mathrm{Im}\:\alpha)=L$, so $\mathcal{Q}$ destabilizes $\mathcal{E}$ and this contradicts the hypothesis.

In the second case we have
\begin{equation*}
\frac{\deg(L)}{r}<\mu(\mathcal{E}).
\end{equation*}
Let $a$ be sufficiently large, $D\in U_a$ and the minimal $\mu$-destabilizing framed quotient
\begin{equation*}
  \begin{tikzpicture}
    \def\x{2}
    \def\y{-1.8}
    \node (A0_0) at (0*\x, 0*\y) {$E\vert_{D}$};
    \node (A0_1) at (1*\x, 0*\y) {$Q_D$};
    \node (A1_0) at (0*\x, 1*\y) {$F\vert_{D}$};
    \path (A0_0) edge [->] node [auto] {$q_D$} (A0_1);
    \path (A0_0) edge [->] node [left=2pt] {$\alpha\vert_{D}$} (A1_0);
  \end{tikzpicture}
\end{equation*}
with $\ker q_D\not\subset \ker\alpha\vert_D.$ By Proposition \ref{prop:quoziente1} (for $\mu$-semistability), $Q_D$ is torsion free, hence there exists an open subscheme $\tilde{D}\subset D$ such that $D\setminus \tilde{D}$ is a closed set of codimension at least two in $D$ and $Q_D\vert_{\tilde{D}}$ is locally free of rank $r.$ Moreover $\ker q_D\vert_{\tilde{D}}\not\subset \ker\alpha\vert_{\tilde{D}}.$ Using the same techniques as in the last part of the proof of \cite[Theorem 7.2.1]{book:huybrechtslehn2010}, we extend $Q_D\vert_{\tilde{D}}$ to a quotient $Q_{\tilde{X}}$ of $\tilde{X}$ which is locally free of rank $r$ with $\det(Q_{\tilde{X}})=L\vert_{\tilde{X}}.$ By construction we have $\ker (E\vert_{\tilde{X}}\rightarrow Q_{\tilde{X}})\not\subset \ker \alpha\vert_{\tilde{X}}$, hence in this way we obtain a quotient $Q$ of $E$ with $\det(Q)=L$ and zero induced framing, such that $\mathcal{Q}$ destabilizes $\mathcal{E}.$
\qedhere\endproof

\subsubsection*{$\mu$-stable case}

Now we want to prove the following generalization of Mehta-Ramana\-than's restriction theorem for $\mu$-stable torsion free sheaves \cite[Theorem 4.3]{art:mehtaramanathan1984}.
\begin{theorem}\label{thm:mr2}
Let $(X,\mathcal{O}_X(1))$ be a nonsingular polarized variety of dimension $d\geq 2.$ Let $F$ be a sheaf on $X$ supported on a divisor $\FD$, which is a locally free $\mathcal{O}_{\FD}$-module. Let $\mathcal{E}=(E,\alpha\colon E\rightarrow F)$ be a $(\FD,F)$-framed sheaf on $X.$ If $\mathcal{E}$ is $\mu$-stable with respect to $\delta_1$, then there is a positive integer $a_0$ such that for all $a\geq a_0$ there is a dense open subset $W_a\subset \vert \mathcal{O}_X(a)\vert$ such that for all $D\in W_a$ the divisor $D$ is smooth, meets transversely the divisor $\FD$ and $\mathcal{E}\vert_D$ is $\mu$-stable with respect to $a\delta_1.$
\end{theorem}
The techniques to prove this theorem are quite similar to the ones used before. By Proposition \ref{prop:extsocle} a $\mu$-semistable $(\FD,F)$-framed sheaf which is simple but not $\mu$-stable has a proper extended framed socle. Thus we first show that the restriction is simple and we use the extended framed socle (rather its quotient) as a replacement for the minimal $\mu$-destabilizing framed quotient.
\begin{proposition}\label{prop:dual}
Let $\mathcal{E}=(E,\alpha)$ be a $\mu$-stable $(\FD,F)$-framed sheaf. For $a\gg 0$ and general $D\in \vert \mathcal{O}_X(a)\vert$ the restriction $\mathcal{E}\vert_D=(E\vert_D,\alpha\vert_D)$ is a simple $(\FD\cap D, F\vert_D)$-framed sheaf on $D.$
\end{proposition} 
To prove this result, we need to define the double dual of a framed sheaf. Let $\mathcal{E}=(E,\alpha)$ be a $(\FD,F)$-framed sheaf; we define a framing $\alpha^{\vee\vee}$ on the double dual of $E$ in the following way: $\alpha^{\vee\vee}$ is the composition of morphisms
\begin{equation*}
E^{\vee\vee}\longrightarrow E^{\vee\vee}\vert_{\FD}\simeq E\vert_{\FD}\stackrel{\alpha\vert_{\FD}}{\longrightarrow} F\vert_{\FD}.
\end{equation*}
Then $\alpha$ is the framing induced on $E$ by $\alpha^{\vee\vee}$ by means of the inclusion morphism $E\hookrightarrow E^{\vee\vee}.$ We denote the framed sheaf $(E^{\vee\vee},\alpha^{\vee\vee})$ by $\mathcal{E}^{\vee\vee}.$ Note that also $\mathcal{E}^{\vee\vee}$ is a $(\FD,F)$-framed sheaf.
\begin{lemma}\label{lem:veestable}
Let $\mathcal{E}=(E,\alpha)$ be a $\mu$-stable $(\FD,F)$-framed sheaf. Then the framed sheaf $\mathcal{E}^{\vee\vee}=(E^{\vee\vee},\alpha^{\vee\vee})$ is $\mu$-stable.
\end{lemma}
\proof
Consider the exact sequence
\begin{equation*}
0\longrightarrow E\longrightarrow E^{\vee\vee}\longrightarrow A\longrightarrow 0
\end{equation*}
where $A$ is a sheaf supported on a closed subset of codimension at least two. Thus $\mathrm{rk}(E^{\vee\vee})=\mathrm{rk}(E)$ and $\deg(E^{\vee\vee})=\deg(E).$ Moreover, since $\alpha=\alpha^{\vee\vee}\vert_{E}$, we have $\mu(\mathcal{E}^{\vee\vee})=\mu(\mathcal{E}).$ Let $G$ be a subsheaf of $E^{\vee\vee}$ and denote by $G'$ its intersection with $E.$ So $\mathrm{rk}(G)=\mathrm{rk}(G')$, $\deg(G)=\deg(G')$ and $\alpha\vert_{G'}=\alpha^{\vee\vee}\vert_{G}.$ Thus we obtain $\mu(G,\alpha^{\vee\vee}\vert_{G})=\mu(G',\alpha\vert_{G'})<\mu(\mathcal{E})=\mu(\mathcal{E}^{\vee\vee}).$
\qedhere\endproof
\begin{definition}
A $d$-dimensional sheaf $G$ on $X$ is \textit{reflexive} if the natural morphism $G\rightarrow G^{\vee\vee}$ is an isomorphism.  
\end{definition}
Recall the following result (see the proof of \cite[Lemma 7.2.9]{book:huybrechtslehn2010}):
\begin{lemma}\label{lem:extrefl}
Let $G$ be a reflexive sheaf. For $a\gg 0$ and $D\in \vert \mathcal{O}_X(a)\vert$ the homomorphism $\mathrm{End}(G)\rightarrow \mathrm{End}(G\vert_D)$ is surjective.
\end{lemma}
\proof[Proof of Proposition \ref{prop:dual}]
For arbitrary $a$ and general $D\in \vert \mathcal{O}_X(a)\vert$ the sheaf $E\vert_D$ is torsion free on $D$ and $E^{\vee\vee}\vert_D$ is reflexive on $D$, moreover the double dual of $E\vert_D$ (as sheaf on $D$) is $E^{\vee\vee}\vert_D$ (cf. \cite[Section 1.1]{book:huybrechtslehn2010}). In addition, $E\vert_D$ is locally free in a neighbourhood of $\FD\cap D$ and the restriction of the framing $\alpha\vert_D$ to $\FD\cap D$ is an isomorphism. We have injective homomorphisms
\begin{eqnarray*}
\delta &\colon& \mathrm{End}(E)\longrightarrow \mathrm{End}(E^{\vee\vee}),\\
\delta_D &\colon& \mathrm{End}(E\vert_D)\longrightarrow \mathrm{End}(E^{\vee\vee}\vert_D).
\end{eqnarray*}
Let $\varphi\in \mathrm{End}(\mathcal{E})$: the image $\varphi^{\vee\vee}=\delta(\varphi)$ of $\varphi$ is an element of $\mathrm{End}(E^{\vee\vee},\alpha^{\vee\vee})$, indeed if $\alpha\circ \varphi=\lambda \alpha$, then we can define an endomorphism of $\mathcal{E}^{\vee\vee}$ in the following way:
\begin{equation*}
  \begin{tikzpicture}
    \def\x{3}
    \def\y{-1.2}
    \node (A0_0) at (0*\x, 0*\y) {$E^{\vee\vee}$};
    \node (A0_1) at (1*\x, 0*\y) {$E^{\vee\vee}$};
    \node (A1_0) at (0*\x, 1*\y) {$E^{\vee\vee}\vert_{\FD}$};
    \node (A1_1) at (1*\x, 1*\y) {$E^{\vee\vee}\vert_{\FD}$};
    \node (A2_0) at (0*\x, 2*\y) {$E\vert_{\FD}$};
    \node (A2_1) at (1*\x, 2*\y) {$E\vert_{\FD}$};
    \node (A3_0) at (0*\x, 3*\y) {$F\vert_{\FD}$};
		\node (A3_1) at (1*\x, 3*\y) {$F\vert_{\FD}$};
    \path (A0_0) edge [->] node [auto] {$\scriptstyle{\varphi^{\vee\vee}}$} (A0_1);
    \path (A0_0) edge [->] node [auto] {$\scriptstyle{}$} (A1_0);
    \path (A0_1) edge [->] node [auto] {$\scriptstyle{}$} (A1_1);
    \path (A1_0) edge [->] node [auto] {$\scriptstyle{\varphi^{\vee\vee}\vert_{\FD}}$} (A1_1);
    \path (A1_0) edge [->] node [auto] {$\scriptstyle{\simeq}$} (A2_0);
    \path (A1_1) edge [->] node [left] {$\scriptstyle{\simeq}$} (A2_1);
    \path (A2_0) edge [->] node [auto] {$\scriptstyle{\varphi\vert_{\FD}}$} (A2_1);
    \path (A2_0) edge [->] node [auto] {$\scriptstyle{\alpha\vert_{\FD}}$} (A3_0);
    \path (A2_1) edge [->] node [left] {$\scriptstyle{\alpha\vert_{\FD}}$} (A3_1);
		\path (A3_0) edge [->] node [auto] {$\scriptstyle{\cdot \lambda}$} (A3_1);
		\path (A0_0) edge [->,bend right=1500] node [left] {$\scriptstyle{\alpha^{\vee\vee}}$} (A3_0);
    \path (A0_1) edge [->,bend left=1500] node [right] {$\scriptstyle{\alpha^{\vee\vee}}$} (A3_1);
  \end{tikzpicture}
\end{equation*}
In the same way it is possible to prove that for $\varphi\in \mathrm{End}(\mathcal{E}\vert_D)$, $\delta_D(\varphi)$ is an element of $\mathrm{End}(\mathcal{E}^{\vee\vee}\vert_D).$ So the restriction of $\delta$ (resp.\ $\delta_D$) to $\mathrm{End}(\mathcal{E})$ (resp.\ $\mathrm{End}(\mathcal{E}\vert_D)$) is still injective. Therefore it suffices to show that $\mathcal{E}^{\vee\vee}\vert_{D}$ is simple for $a\gg 0$ and general $D.$ By Lemma \ref{lem:veestable}, $\mathcal{E}^{\vee\vee}$ is $\mu$-stable, hence by point (c) of Corollary \ref{cor:semist} it is simple. By Lemma \ref{lem:extrefl}, the homomorphism $\chi\colon \mathrm{End}(E^{\vee\vee})\rightarrow \mathrm{End}(E^{\vee\vee}\vert_D)$ is surjective for $a\gg 0$ and general $D.$ Since for $\varphi\in \mathrm{End}(\mathcal{E}^{\vee\vee})$, $\chi(\varphi)$ is an element of $\mathrm{End}(\mathcal{E}^{\vee\vee}\vert_D)$, we have that the map
\begin{equation*}
\chi\vert_{\mathrm{End}(\mathcal{E}^{\vee\vee})}\colon \mathrm{End}(\mathcal{E}^{\vee\vee})\rightarrow \mathrm{End}(\mathcal{E}^{\vee\vee}\vert_D)
\end{equation*} 
is also surjective. Thus $ \mathrm{End}(\mathcal{E}\vert_D)=\mathrm{End}(\mathcal{E}^{\vee\vee}\vert_D)\simeq k.$
\qedhere\endproof
\begin{remark}\label{rem:notdest}
Since $\mathcal{E}$ is $\mu$-stable with respect to $\delta_1$, we have $\deg(\mathrm{Im}\:\alpha)> \delta_1.$ For any positive integer $a$ and $D\in \Pi_a$, $\deg(\mathrm{Im}\:\alpha\vert_D)=a\deg(\mathrm{Im}\:\alpha)> a\delta_1$, hence $\ker\alpha\vert_D$ is not framed $\mu$-destabilizing. Therefore it cannot be the extended framed socle of $\mathcal{E}\vert_D.$ Thus, for any positive integer $a$ and general $D\in \Pi_a$, the quotient of $\mathcal{E}\vert_D$ by its extended framed socle has positive rank.\triend
\end{remark}
Now we shall prove Theorem \ref{thm:mr2}. Since the arguments we shall use are similar to that of the nonframed case, we give only a sketch of the proof and we refer to the proof of \cite[Theorem 7.2.8]{book:huybrechtslehn2010}.
\proof[Proof of Theorem \ref{thm:mr2}]
Let $a_0\geq 3$ be an integer such that for all $a\geq a_0$ and a general $D\in \Pi_a$, the restriction $\mathcal{E}\vert_D$ is $\mu$-semistable with respect to $a\delta_1$ and simple (cf. Proposition \ref{prop:dual}). 

Assume that the theorem is false: for all $a\geq a_0$ and general $D\in \Pi_a$, $\mathcal{E}\vert_D$ is not $\mu$-stable with respect to $a\delta_1.$ Then $\mathcal{E}\vert_{D_\eta}$ is not geometrically $\mu$-stable for the divisor $D_{\eta}$ associated to the generic point $\eta\in \Pi_a$, i.e., the pull-back to some extension of $k(\eta)$ is not $\mu$-stable (cf. Corollary \ref{cor:geostable}). Hence $\mathcal{E}\vert_{D_{\eta}}$ is not $\mu$-stable. Since $\mathcal{E}\vert_{D_{\eta}}$ is simple, by Proposition \ref{prop:extsocle} the extended socle of $\mathcal{E}\vert_{D_{\eta}}$ is a proper $\mu$-destabilizing framed subsheaf (different from $\ker\alpha\vert_{D_\eta}$ by Remark \ref{rem:notdest}). Consider the corresponding positive rank quotient sheaf $Q_\eta$, with induced framing $\beta_\eta$: we can extend it to a quotient $q^*E\rightarrow Q_a$ over all of $Z_a.$ 

Let $W_a$ be the dense open subset of points $D\in \Pi_a$ such that 
\begin{itemize}
 \item $D$ is smooth, meets transversely the divisor $\FD$ and $E\vert_D$ is torsion free, 
\item $Q_a$ is flat over $W_a$ and $\epsilon\left((\tilde{\alpha}_a)\vert_D\right)=\epsilon(\beta_\eta)$, where we denote by $\tilde{\alpha}_a$ the induced framing on $Q_a.$
\end{itemize}
Thus $Q_a\vert_D$ is a sheaf of positive rank such that, with the induced framing, it is a $\mu$-destabilizing framed quotient for all $D\in W_a.$

Using the same arguments of proof of \cite[Theorem 7.2.8]{book:huybrechtslehn2010}, one can prove that there is a line bundle $L$ on $X$ and an integer $0<r<\mathrm{rk}(E)$ such that for $a\gg 0$ and for general $D\in W_a$
\begin{equation*}
 \mu(Q_a\vert_D,\tilde{\alpha}_a\vert_D)=\frac{\deg(L\vert_D)-\epsilon(a)a\delta_1}{r}=a\Big(\frac{\deg(L)-\epsilon(a)\delta_1}{r}\Big)=\mu(E\vert_D,\alpha\vert_D)=a\mu(E,\alpha),
\end{equation*}
hence
\begin{equation*}
\frac{\deg(L)-\epsilon(a)\delta_1}{r}=\mu(E,\alpha).
\end{equation*}
Using the arguments at the end of the proof of the restriction theorem for $\mu$-semistable framed sheaves, one can show that this suffices to construct a $\mu$-destabilizing framed quotient $\mathcal{E}\rightarrow \mathcal{Q}$ for sufficiently large $a.$ This contradicts the assumptions of the theorem.
\qedhere\endproof

\end{document}